\documentclass[review]{elsarticle}

\usepackage{lineno,hyperref}
\usepackage{subfigure}
\usepackage{amsmath, amssymb}
\usepackage{amsbsy}
\usepackage{mathptmx}
\usepackage{bm}
\usepackage{makecell} 

\journal{Elsevier}

\begin{document}

\begin{frontmatter}

\title{Modified multi-dimensional limiting process with enhanced shock stability on unstructured grids}

\author[mymainaddress,mysecondaryaddress]{Fan Zhang\corref{mycorrespondingauthor}}
\ead{a04051127@mail.dlut.edu.cn}
\cortext[mycorrespondingauthor]{Corresponding author.}

\author[mymainaddress,mysecondaryaddress]{Jun Liu}

\author[mysecondaryaddress]{Biaosong Chen}

\address[mymainaddress]{State Key Laboratory of Aerodynamics, Mianyang, Sichuan, China }
\address[mysecondaryaddress]{State Key Laboratory of Structural Analysis for Industrial Equipment, Dalian University of Technology, Dalian, China}

\begin{abstract}
The basic concept of multi-dimensional limiting process (MLP) on unstructured grids is inherited and modified for improving shock stabilities and reducing
numerical dissipation on smooth regions. A relaxed version of MLP condition, simply named as weak-MLP, is proposed for reducing dissipation.
Moreover, a stricter condition, that is the strict-MLP condition, is proposed to enhance the numerical stability. The maximum/minimum principle is fulfilled by both the strict- and
weak-MLP condition. A differentiable pressure weight function is applied for the combination of two novel conditions, and thus the modified limiter is named
as MLP-pw(pressure-weighted). A series of numerical test cases show that MLP-pw limiter has improved stability and convergence,
 especially in hypersonic simulations. Furthermore, the limiter also shows lower dissipation in regions without significant pressure transition. Therefore,
 MLP-pw limiter can capture contact discontinuity and expansion accurately.
\end{abstract}


\end{frontmatter}
 
\section{Introduction} \label{sec:Intro}

Unstructured grids are commonly used for spatial discretizations by current industrial computational fluid dynamics codes
that simulate aerodynamics or gas dynamics phenomena. The advantages of unstructured grids include the conveniences in automatic grid generation \cite{Bowyer1981,Watson1981,Zhang2013},
grid adaptation \cite{Drikakis1997,Alauzet2007,Harris2014}, moving mesh techniques \cite{Luke2012,LIU2012}, for complex geometries and flow phenomena.
The application of unstructured grids significantly facilitates the aforementioned aspects.
However, the accuracy and stability of unstructured schemes are usually challenged by the irregularity of grid
connectivity and deterioration of grid quality \cite{AndreasMack2002,Zhang2016a}, which are inevitable for the automatic discretization of complicate geometries.
Especially, the simulations for transonic and supersonic flows require good approximation of nonlinear multi-dimensional
physical phenomena such as shock waves, shock waves interaction, and shock-vortex interaction, and thus good accuracy and stability are indispensable.

As a key factor that affects spatial accuracy and stability, slope limiter, or for short, limiter, has been investigated for decades.
As well known, second-order or higher than second-order schemes suffer from numerical oscillations across discontinuities, a typical one of which is shock wave  \cite{Anderson1974,Beam1976}. Therefore, limiter is used to suppress these oscillations while keeping second-order accurate reconstruction in smooth regions of flow fields.
On structured grids, the finite difference method (FDM) and finite volume method (FVM) have been applied along with mature limiting method based on
solid theories. The typical strategy is MUSCL (Monotonic Upstream-Centered Scheme for Conservation Laws) scheme \cite{Leer1979} with limiters that subject to TVD (Total variation diminishing) condition \cite{Harten1983,Harten1984,Sweby1984}.
However, these structured schemes can not be extended to unstructured grids directly due to various reasons.
Firstly, the schemes for structured grids are usually developed based on one-dimensional analysis and extended to multi-dimensional structured grids by dimensional-slitting, which is infeasible for unstructured grids. Secondly, one-dimensional principles, for instance, the TVD condition, are not necessary
feasible in multi-dimensional unstructured grids. A counterexample of Jameson had shown that a flow field on which the total variation is smaller could be more
oscillatory than a flow field on which the total variation is larger \cite{Jameson1995}. Furthermore, a scheme in TVD condition will causes
accuracy deterioration in extrema even in smooth region, and thus the TVB \cite{Shu1987mc} and ENO \cite{Harten1987} schemes were developed.

By extending Spekreijse's monotone condition \cite{Spekreijse1987mc}, Barth and Jespersen designed a limiter on unstructured grids \cite{Barth1989},
which modifies the piecewise linear distribution at each control volume. Barth-Jespersen limiter removes local extrema and insures stability. However,
this limiter shows similar effects as that of TVD condition, which reduces accuracy on smooth extrema. Furthermore, the limit function
of Barth and Jespersen is non-differentiable, and thus the convergence is less satisfactory. Therefore, an improvement was introduced by
Venkatakrishnan \cite{VENKATAKRISHNAN1993}, who used a differentiable function similar to van Albada limiter \cite{vanAlbada1982} which is designed on structured grids. Venkatakrishnan limiter archives better convergence compared with Barth-Jespersen limiter. Whereas, Venkatakrishnan limiter isn't strictly monotone, and thus it may produce oscillations across shock waves. Generally speaking, Barth-Jespersen limiter and Venkatakrishnan limiter have been
commonly applied on unstructured grids since their inventions.

Many researches have been focusing on the improvement of limiters. In order to reduce the dissipations of two aforementioned unstructured limiters,
a strategy was introduced, which is turn off limiter in subsonic regions. Nejat and Ollivier-Gooch introduced hyperbolic tangent function
in their application of Venkatakrishnan limiter, by which the limiter only activates in limited regions \cite{Nejat2008}.
Michalak and Ollivier-Gooch further improved this method \cite{Michalak2009}.
Thereafter, Kitamura and Shima introduced the concept of second limiter, which also uses a hyperbolic tangent function to turn off limiter in stagnation or subsonic regions, but removes predefined parameters \cite{Kitamura2012a}. It was proved by numerical results that second limiters can reduce dissipations effectively.

A relatively new method on unstructured grids is MLP (Multi-dimensional Limiting Process) limiter, which was first introduced on structured grids \cite{Kim2005a,Yoon2008}. By using the MLP condition which satisfies maximum/minimum principles, MLP limiter properly introduces multidimensional information. Therefore, the method has been showing better accuracy, robustness and convergence in various circumstances.
Park, et al. designed unstructured MLP limter \cite{Park2010}.
Thereby, Park and Kim \cite{Park2012} had constructed three-dimensional unstructured MLP limiter and proved that the limiter
obeys LED (Local Extremum Diminishing) condition \cite{Jameson1995}. Gerlinger designed a low dissipation MLP limiter, MLP$^{ld}$,
on structured grids, and simulated combustion problem \cite{Gerlinger2012}.
Do, et al. defined a low dissipation MLP limiter for central-upwind Schemes \cite{Do2016}.
Kang, et al. \cite{Kang2010} reduced dissipation by only turning on the MLP limiter in the vicinity of shock waves/nonlinear discontinuity.
MLP limiter had also been developed for higher order unstructured numerical schemes \cite{Park2014,Park2016}.
Li, et al. developed a multi-dimensional limiter, WBAP, which modifies the gradients by a component by component approach \cite{Li2011,Li2012}.
This method is not rotationally invariant but shows good accuracy, robustness and convergence performance in numerical tests.

In spite of the successful applications of former researchers and the authors, there is still room for MLP limiter to improve the stability and convergence,
especially for hypersonic flow simulations. Therefore, the presented research is focusing on this topic. This paper is organized as follows.
The finite volume method and spatial reconstruction are briefly described in Section \ref{sec:FVM}.
Then, the Barth-Jespersen limiter, Venkatakrishnan limiter and MLP limiter are introduced in section \ref{sec:limiters}, where the differences are
emphasised. In Section \ref{sec:Imp}, the presented modifications on MLP limiter are formulated.
A series of numerical test cases along with corresponding discussions are given in section \ref{sec:tests}. Finally, section \ref{sec:Conclusions} concludes the whole work.

\section{Finite volume method and second-order reconstruction} \label{sec:FVM}

The discretization for the compressible Navier-Stokes equations is introduced as follows. The integral form of the equations is
\begin{equation} \label{eq:GoverEq}
\int\limits_\Omega\frac{\partial\mathbf{Q}}{\partial t}\mathrm{d}\Omega+\int\limits_{\partial \Omega}[\mathbf{F}_c(\mathbf{Q})-\mathbf{F}_v(\mathbf{Q})]\cdot\mathbf{n}\mathrm{d}S=\mathbf{0},
\end{equation}

\noindent where $\mathbf{Q}$ are the conservative variables in the flow field, $\mathbf{F}_c(\mathbf{Q})$ is convective flux, and $\mathbf{F}_v(\mathbf{Q})$ is viscous flux, which could be solved
 by using a central scheme for unstructured grids \cite{Blazek2005}.
In this paper, solutions of the convective flux are emphatically investigated. Therefore, in the following discussions the $\mathbf{F}_v(\mathbf{Q})$ term is neglected, and thus the equations are simplified as Euler equations. $\mathbf{Q}$ and $\mathbf{F}_c(\mathbf{Q})$ are given as

\begin{equation} \label{eq:Var_Flux}
\mathbf{Q} =
 \begin{bmatrix}
\rho   \\
\rho u \\
\rho v \\
\rho E
\end{bmatrix},\\
\mathbf{F}_c(\mathbf{Q})=
 \begin{bmatrix}
\rho V_n  \\
\rho uV_n +pn_x\\
\rho vV_n +pn_y \\
 \rho HV_n
\end{bmatrix},
\end{equation}

\noindent where $V_n=\mathbf{V} \cdot \mathbf{n}=(un_x+vn_y)$. $E$ is the total energy, $H$ is the enthalpy, given as
\begin{equation}
E=\frac{1}{\gamma-1}\frac{p}{\rho}+\frac{1}{2}(u^2+v^2),
\end{equation}
\begin{equation}
H=E+\frac{p}{\rho},
\end{equation}
\noindent where $\gamma$ is the ratio of specific heat. For air at moderate pressures and temperatures one uses $\gamma=1.4$.

The governing equations are discretized by using cell-centered finite volume formulation which is applied to a polygon computational cell $i$
sharing a interface $k$ with a neighbouring cell $j$.
Therefore, the spatial discretization at cell $i$ for the Euler equations can be expressed as
\begin{equation} \label{eq:SpatialDiscrete}
\frac{\partial}{\partial t}(\mathbf{Q}\Omega)_i=-(\sum_{k=1}^{N_f}\mathbf{F}_{c, k}\cdot\mathbf{n}_k S_k)_i,
\end{equation}

\noindent where $S_k=|\partial\Omega_k|$ is the interface area, $\mathbf{n}_k$ is the unit norm vector outward from the interface, $N_f$ is the interface number of cell $i$. Although the exact convective flux function $\mathbf{F}_{c, k}$ is nonlinear, it is usually
solved by a linearized numerical flux instead of the exact formula \cite{godunov1959}.
Furthermore, the numerical flux function could be simplified as an one-dimensional scheme that calculates in the direction of vector $\mathbf{n}_k$.
In fact, upwind schemes, such as FDS (Flux Difference Splitting) scheme or FVS (Flux Vector Splitting) scheme, are mostly designed based on one-dimensional hypothesis.
The FDS scheme or FVS scheme could be defined as a function of conservative variables $\mathbf{Q}$, and thus the flux is given as
\begin{equation} \label{eq:Flux}
\mathbf{F}_{c,k}=\mathbf{F}_\text{FDS/FVS}(\mathbf{Q}_{k}^+, \mathbf{Q}_{k}^-,\mathbf{n}_k),
\end{equation}

\noindent where the superscript $(\cdot)^{\pm}$ denote the left and right values of interface $k$ respectively.
In the following paragraphs, the subscripts $c$ and $k$ are neglected for simplicity.

The cell interface values are extrapolated from the cell centre values by using gradient $\nabla q$:
\begin{equation} \label{eq:Second_order}
\begin{aligned}
q_{k}^+ &= {q}_i+ \phi_i\nabla q_i\cdot\Delta\mathbf{r}_{ik},  \\
q_{k}^- &= {q}_j+ \phi_j\nabla q_j\cdot\Delta\mathbf{r}_{jk},
\end{aligned}
\end{equation}

\noindent where $\Delta(\cdot)_{ik}=(\cdot)_k-(\cdot)_i$ and $q$ could be any of the conservative variables. $\nabla q$ is calculated by nodal averaging procedure
\cite{Frink1991} and
 Gauss-Green scheme \cite{Barth1989}, and
the slope limiter value $\phi$ is employed to suppress oscillations at captured discontinuities. In the following sections, the calculation of $\phi$ will be investigated.
Reconstruction becomes conservative if the integration of $q$ over a cell equals to the cell-averaged value
\begin{equation}
\overline{q}=\frac{1}{|\Omega|}\int\limits_\Omega q\mathrm{d}\Omega.
\end{equation}

The time derivative in Eq.\ref{eq:SpatialDiscrete} could be solved by many explicit and implicit schemes. Due to the limited topic of the presented article, temporal
solutions will not be further discussed.

\section{Limiters} \label{sec:limiters}

Barth-Jespersen limiter, Venkatakrishnan limiter and MLP limiter are briefly introduced in the section. The functions of Barth-Jespersen and Venkatakrishnan could be utilised by MLP limiter and the modified limiter of the presented article.

\subsection{Barth-Jespersen and Venkatakrishnan Limiters}

Barth-Jespersen limiter and Venkatakrishnan limiter are two typical and common used limiters on unstructured grids.
In fact, these two limiter are following similar formulation.

For the calculation of $\phi_i$ that limits a scalar variable in a cell $i$, the Barth-Jespersen limiter is given as

\begin{equation} \label{eq:BJ}
\begin{aligned}
\phi_\text{BJ} =\min \left\{
\begin{aligned}
&f_\text{BJ} \left(\frac{q^{\text{max}}_i-q_i}{q_\text{t}-q_i} \right), &\text{if} \quad q_\text{t}-q_i > 0 \\
&f_\text{BJ} \left(\frac{q^{\text{min}}_i-q_i}{q_\text{t}-q_i} \right), &\text{if} \quad q_\text{t}-q_i < 0 \\
&1, &\text{if} \quad q_\text{t}-q_i = 0
\end{aligned}
\right.
\end{aligned}
\end{equation}
\noindent where the subscript \emph{t} indicates a \emph{test} value, which could be different based on definitions, and the function $f_\text{BJ}$ is defined as
\begin{equation} \label{eq:BJ_f}
f_\text{BJ}\left(\frac{\Delta_+}{\Delta_-}\right)=\min\left(1, \frac{\Delta_+}{\Delta_-}\right).
\end{equation}

\noindent where the $\Delta_+=q^{\text{max,min}}-q$ and $\Delta_-=q_\text{t}-q$.
This minimum function will limit the gradient, i.e. $\phi_i<1$, if the \emph{test} value is larger than the maximum, $q^{\max}_i$, or smaller than the minimum, $q^{\min}_i$.
It had been commonly reported that the non-differentiable function $f_\text{BJ}$ is a major drawback that may cause accuracy lost in smooth regions and convergence
deterioration in steady state computations.

Therefore, an improvement, the limiter of Venkatakrishnan, is given as
\begin{equation} \label{eq:Venkatakrishnan}
\begin{aligned}
\phi_\text{V} = \min\left\{
\begin{aligned}
&f_\text{V} \left(\frac{q^{\text{max}}_i-q_i}{q_\text{t}-q_i} \right), &\text{if} \quad q_\text{t}-q_i > 0 \\
&f_\text{V}  \left(\frac{q^{\text{min}}_i-q_i}{q_\text{t}-q_i} \right), &\text{if} \quad q_\text{t}-q_i < 0 \\
&1, &\text{if} \quad q_\text{t}-q_i = 0
\end{aligned}
\right.
\end{aligned}
\end{equation}
\noindent where the function $f_\text{V}$ is
\begin{equation} \label{eq:Venkata_f}
f_\text{V}\left(\frac{\Delta_+}{\Delta_-}\right)=\frac{1}{\Delta_{-}}\left[\frac{({\Delta_+}^2+\varepsilon^2)\Delta_-+2{\Delta_-}^2\Delta_+}
 {{\Delta_+}^2+2{\Delta_-}^2+\Delta_+\Delta_-+\varepsilon^2}\right],
\end{equation}
\noindent which is differentiable. The small parameter $\varepsilon^2$ in limit function is defined as
\begin{equation}
\varepsilon^2=(K\Delta h)^3,
\end{equation}
\noindent by which the limitation effect of limiter is tunable. Parameter $\Delta h$ is the cell scale, and $K$ usually is defined by user to adjust numerical dissipation. Parameter $K$ is used to tune the restriction of the limiter and was evaluated in \cite{VENKATAKRISHNAN1993}. In general, if $K=0$, the limiter will be very dissipative because it's activated even in smooth region. Conversely, the limiter will be actually turned off if $K\gg 1$.

The difference between Barth-Jespersen and Venkatakrishnan limiter is produced by introducing the function $f_\text{V}$, by which the convergence is improved,
and the dissipation is reduced, but the monotonicity will not be strictly guaranteed \cite{VENKATAKRISHNAN1993}.

In Barth-Jespersen or Venkatakrishnan limiters, the maximum and minimum values among the direct neighbouring cells are given as
\begin{equation}
q^{\text{max}}_i=\max(q_i, \max\limits_{j\in V(i)}q_j ), \quad
q^{\text{min}}_i=\min(q_i, \min\limits_{j\in V(i)}q_j  ),
\end{equation}
\noindent where the subscript $V(i)$ indicates a set of cells connected with cell $i$ by a common interface, and the $V$ means \emph{Volume}.

In order to avoid any new extrema in each cell, the following condition is expected
\begin{equation} \label{eq:max_min_1}
q^{\text{min}}_i\le q_{\text{t}} \le q^{\text{max}}_i,
\end{equation}
\noindent where the $q_{\text{t}}$ could be any reconstructed variables within the cell. Two different choices could be used for the definition of $q_{\text{t}}$, in a discretized manner. One may chooses to restrict the variable distribution in the whole cell by applying Eq.\ref{eq:max_min_1} at each \emph{vertex}, and thus the following formula could be given
\begin{equation} \label{eq:Limit_Node}
q_\text{t}^{(v)}=q_i+\nabla q_i \cdot \Delta\mathbf{r}_{il}, \quad l\in v(i),
\end{equation}
\noindent where the lower-case \emph{v} indicates \emph{vertex}, and thus the superscript $(v)$ means the value is defined at a vertex, and $v(i)$ indicates the set of grid vertexes of cell $i$. If we only expect that the values at interface centres satisfy the condition of Eq.\ref{eq:max_min_1}, the following equation will be given
\begin{equation} \label{eq:Limit_Face}
q_\text{t}^{(f)}=q_i+\nabla q\cdot \Delta\mathbf{r}_{ik}, \quad \partial\Omega_k\subset \partial\Omega_i,
\end{equation}
\noindent where the lower-case superscript $(f)$ means the value is defined at an interface centre. It is obvious that the Eq.\ref{eq:Limit_Node} will be more strict and dissipative because a linear reconstruction always
shows maximum and minimum at the vertexes which are the most distant points from cell centre. In the numerical cases of this article,
the definition in Eq.\ref{eq:Limit_Node} will be applied for Barth-Jespersen limiter and Venkatakrishnan limiter, of which the stability will be expected
to be enhanced.

\subsection{Unstructured MLP Limiter}

The MLP limiter on unstructured grids is somehow simpler than its structured version. It should be noted that the MLP limiter showed different standpoint
compared with that of Barth-Jespersen's or Venkatakrishan's method. A simple but important condition is
\begin{equation} \label{eq:MLP_1}
q_{V(l)}^\text{min} \le q_{l} \le q_{V(l)}^\text{max}, \quad l\in v(i),
\end{equation}
\noindent where the subscript $l$ indicates a \emph{vertex} which belongs to the vertexes set $v(i)$ of cell $i$, and the $V(l)$ indicates the set of all the cells connected to vertex $l$. And
\begin{equation} \label{eq:MLP_1_1}
q_{V(l)}^\text{max} = \max\limits_{j\in V(l)}(q_j), \quad
q_{V(l)}^\text{min} = \min\limits_{j\in V(l)}(q_j).
\end{equation}

\begin{figure}
\begin{center}
\includegraphics[width=5cm,angle=-90]{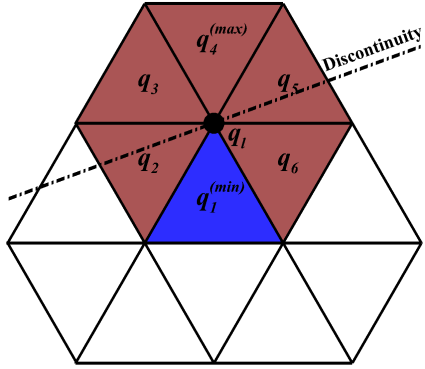}
\caption{Stencils involved in MLP condition.}
\label{fig:f:MLP_limit}
\end{center}
\end{figure}

Eq.\ref{eq:MLP_1} is the so call MLP condition. For each vertex of a cell, all the cell-averaged values sharing this vertex are utilised for detecting flow phenomena, including discontinuities,
 as shown in Fig.\ref{fig:f:MLP_limit}.

By applying the MLP condition to Eq.\ref{eq:Second_order}, the permissible range of $\phi$ is calculated by
\begin{equation} \label{eq:MLP_2}
\frac{q_{V(l)}^\text{min}-q_i}{\nabla q \cdot\Delta\mathbf{r}_{il}}  \le \phi_{l} \le \frac{q_{V(l)}^\text{max}-q_i}{\nabla q \cdot\Delta\mathbf{r}_{il}}.
\end{equation}
\noindent The minimum $\phi_i=\min\limits_{l\in v(i)}(\phi_{l})$ of a cell will be applied for linear reconstruction. The following maximum/minimum principle
\begin{equation} \label{eq:MLP_3}
q_{\mathcal{V}(i)}^{\text{min},n} \le q_{i}^{n+1} \le q_{\mathcal{V}(i)}^{\text{max},n}
\end{equation}
\noindent is satisfied by the linear reconstruction, where
\begin{equation}
q_{\mathcal{V}(i)}^{\max} = \max\limits_{l\in v(i)}(q_{V(l)}^\text{max}), \quad
q_{\mathcal{V}(i)}^{\min} = \min\limits_{l\in v(i)}( q_{V(l)}^\text{min}).
\end{equation}
\noindent  MLP limiter will be utilising a wider range of flow information, that is the set of common vertex neighbouring cells $j\in\mathcal{V}(i)$, compared with Barth-Jespersen and Venkatakrishnan limiters, which use only the direct neighbouring cells, $j\in {V}(i)$.
Because more information will be applied by the limiter, the limiter will be more accurate and less sensitive
 to grid perturbation. For detailed discussions one may refer to article \cite{Kim2005a} and \cite{Park2010}.

Eventually, the formula of MLP limiter is
\begin{equation} \label{eq:MLP_4}
\begin{aligned}
\phi_\text{MLP} =\min\limits_{l\in v(i)} \left\{
\begin{aligned}
&f_\text{MLP} \left(\frac{q_{V(l)}^\text{max}-q_i}{q_{l}-q_i} \right), &\text{if} \quad q_{l}-q_i > 0 \\
&f_\text{MLP} \left(\frac{q_{V(l)}^\text{min}-q_i}{q_{l}-q_i} \right), &\text{if} \quad q_{l}-q_i < 0 \\
&1, &\text{if} \quad q_{l}-q_i = 0
\end{aligned}
\right.
\end{aligned}
\end{equation}
\noindent where the $q_{l}$ is a vertex value calculated by unlimited linear reconstruction from the centre of cell $i$.    The function $f_\text{MLP}$ could be the function of Barth-Jespersen limiter, $f_\text{BJ}$, or that of Venkatakrishnan Limiter, $f_\text{V}$. By using these two different functions, MLP limiter could be showing different performance \cite{Park2010,Park2012}.
In the numerical cases of this article, due to its better performance in convergence, the function $f_{\text{V}}$ will be applied for MLP limiter.

\section{A modification: MLP-pw limiter} \label{sec:Imp}

Although modifications are made in the presented scheme, the basic idea is following that of MLP limiter.
The stencils in Fig.\ref{fig:f:MLP_limit} are also implemented in the modified limiter, which is considered as a development of MLP.
The differences are produced because of the purpose of improving shocks stability and reducing dissipation on continuous region or contact discontinuity.

\subsection{A relaxed version of MLP condition} \label{sec:relex}
Several researches had managed to reduce the dissipation of slope limiters, on both structured and unstructured grids.
A typical and effective example on unstructured grids is the second limiter of Kitamura and Shima \cite{Kitamura2012a}, which turns off the limiter in subsonic regions.
However, on unstructured grids, the slope limiter is securing the numerical schemes from unphysical spatial reconstruction, which is not only the oscillations
across discontinuity, but also could be geometrical monotonicity violation \cite{Shima2013a}. Therefore, an alternative strategy will be developed for reducing dissipation.

The strategy used here is to perform limitation at each interface, and thus the following formula is presented
\begin{equation} \label{eq:f_extr}
\overline{q}_{k}^{\text{max}}=  \left. \sum\limits_{l\in v(k)} q_{V(l)}^{\text{max}}  \middle /   n_{v(k)} \right., \quad
\overline{q}_{k}^{\text{min}}= \left.  \sum\limits_{l\in v(k)} q_{V(l)}^{\text{min}}   \middle / n_{v(k)}\right. ,
\end{equation}
\noindent where the $q_{V(l)}^{\max,\min}$ were defined in Eq.\ref{eq:MLP_1_1}, $n_{v(k)}$ is the number of elements in the set $v(k)$, which includes all the vertexes of interface $k$. For the two-dimensional cases, $n_{v(k)}\equiv 2$ because the interface (line) only connects two vertexes. Then the condition is given as
\begin{equation} \label{eq:wMLP}
\overline{q}_{k}^{\text{min}}\le q_k^\pm \le \overline{q}_{k}^{\text{max}},
\end{equation}
\noindent where $q_k^\pm $ is reconstructed values at each side of interface $k$. This condition is named as weak-MLP condition because it is less restrictive compared with the (original) MLP condition.

It is necessary to prove the numerical dissipation of weak-MLP condition is indeed reduced. Therefore, the following lemma for linear reconstruction is given
based on two-dimensional assumption.

\noindent \textbf{Lemma.} A linear reconstruction that satisfies the weak-MLP condition in Eq.\ref{eq:wMLP} is not more diffusive, and could be less diffusive, compared with the reconstruction satisfying MLP condition.

\begin{figure}
\begin{center}
\includegraphics[width=5cm,angle=-90]{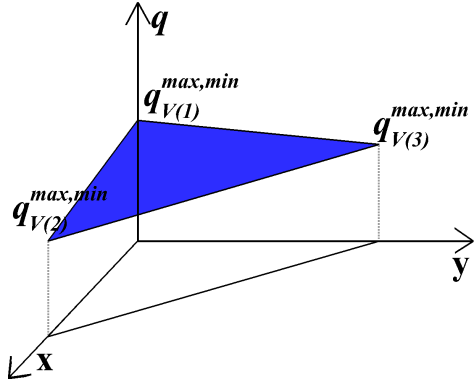}
\caption{Function plane constructed by $q_{V(l)}^{\max,\min}$.}
\label{fig:f:q_face}
\end{center}
\end{figure}

\noindent \textbf{Proof.}  Assume the reconstruction of cell $i$ is made on a triangular cell, and $v(i)=\{1,2,3\}$.
At each vertex $l\in v(i)$, there are a maximum, $q_{V(l)}^{\max}$, and a minimum, $q_{V(l)}^{\min}$, which are already given. Therefore, an unique linear
distribution could be defined for all the maximums or all the minimums,
 \begin{equation} \label{eq:plane_m_m}
q^{(\max)}(\mathbf{r})=\overline{q}^{\max}_i+\nabla q^{\max}(\mathbf{r}-\mathbf{r}_{i}), \quad
q^{(\min)}(\mathbf{r})=\overline{q}^{\min}_i+\nabla q^{\min}(\mathbf{r}-\mathbf{r}_{i}),
\end{equation}
\noindent where the $\overline{q}_i^{\max,\min}$ are the averages of maximums and minimums respectively. $q^{(\max,\min)}(\mathbf{r})$ are linear functions,
which means that their distribution are planes in function space, as shown in Fig.\ref{fig:f:q_face}.

 Based on the definition of $q_{V(l)}^{\max,\min}$ in Eq.\ref{eq:MLP_1_1}, the cell centre value $q_i$ will be satisfying
 \begin{equation} \label{eq:plane_m_m2}
 \max\limits_{l\in v(i)}(q_{V(l)}^{\min}) \le q_i \le \min\limits_{l\in v(i)}(q_{V(l)}^{\max}).
\end{equation}
\noindent If $\max\limits_{l\in v(i)}(q_{V(l)}^{\max})=\min\limits_{l\in v(i)}(q_{V(l)}^{\max})$, two possibilities could be given:

\noindent p1. $q_i$ is the maximum among all the common vertex neighbouring cells, $\mathcal{V}(i)$,

\noindent p2. cell centre values in neighbouring cells are all larger than $q_i$ and equal to each other.

\noindent Therefore, one of the following equation will be valid
 \begin{equation}
 q_i=q^{(\max)}(\mathbf{r}_i)=\overline{q}^{\max}_i, \quad  q_i=q^{(\min)}(\mathbf{r}_i)=\overline{q}^{\min}_i.
\end{equation}

\noindent If $\max\limits_{l\in v(i)}(q_{V(l)}^{\min})=\min\limits_{l\in v(i)}(q_{V(l)}^{\min})$, similar result will be given.
In such circumstances, MLP condition and weak-MLP condition are showing not difference.

Otherwise, $q_i$ will be satisfying the following inequality
 \begin{equation} \label{eq:plane_m_m3}
 q^{(\min)}(\mathbf{r}_i)< q_i < q^{(\max)}(\mathbf{r}_i).
\end{equation}

Considering the maximum constraint $q^{(\max)}(\mathbf{r})$, the following deduction could be given. A limit value $\phi_\text{MLP}$ is given for a linear reconstruction $\phi \nabla q$, which will be a plane in function space. Three circumstances that satisfy MLP condition are listed follow:

\begin{figure}
\begin{center}
\includegraphics[width=4cm,angle=-90]{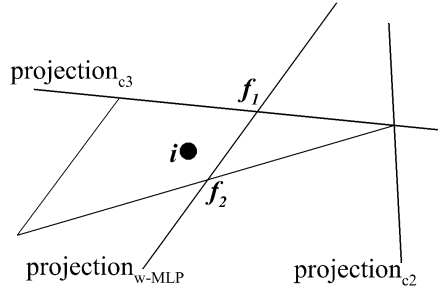}
\caption{Projections on $x-y$ space of intersection lines in function space.}
\label{fig:f:proj_intersection}
\end{center}
\end{figure}

\noindent c1. If the reconstruction plane is not cutting the plane of $q^{(\max)}(\mathbf{r})$ within the cell, it will be not limitation taking place.

\noindent c2. If the reconstruction plane is cutting the plane of $q^{(\max)}(\mathbf{r})$, and the projection of intersection line is
across and only across the cell at one of the vertexes of the cell (projection$_\text{c2}$ in Fig.\ref{fig:f:proj_intersection}), the MLP condition will not be actived except at this vertex.
Therefore, except at the vertex, following inequality is valid
 \begin{equation} \label{eq:plane_m_m4}
q_\text{MLP}(\mathbf{r})=q_i+\phi_\text{MLP} \nabla q\cdot (\mathbf{r}-\mathbf{r}_i)<q^{(\max)}(\mathbf{r}).
\end{equation}
\noindent The weak-MLP condition will be satisfied as well.

\noindent c3. If the reconstruction plane is cutting the plane of $q^{(\max)}(\mathbf{r})$, and the projection of intersection line coincide with one of the interface of the cell (projection$_\text{c3}$ in Fig.\ref{fig:f:proj_intersection}), the MLP condition will be actived at the two vertexes of the corresponding interface. Except at the intersection line/interface, Eq.\ref{eq:plane_m_m4} is satisfied,
and thus the weak-MLP condition is satisfied as well.

Then, a linear reconstruction that satisfies weak-MLP is given as
 \begin{equation} \label{eq:plane_m_m5}
q_\text{w-MLP}(\mathbf{r})=q_i+\phi_\text{w-MLP} \nabla q\cdot (\mathbf{r}-\mathbf{r}_i).
\end{equation}
\noindent The reconstructed plane in function space is cutting the plane of $q^{(\max)}(\mathbf{r})$ in a line of which the projection is across the centres of
 interface $f_1$ and $f_2$ (projection$_\text{w-MLP}$ in Fig.\ref{fig:f:proj_intersection}), and then
 \begin{equation} \label{eq:plane_m_m6}
q_\text{w-MLP}(\mathbf{r}_{f_1})=q^{(\max)}(\mathbf{r}_{f_1}), \quad q_\text{w-MLP}(\mathbf{r}_{f_2})=q^{(\max)}(\mathbf{r}_{f_2}),
\end{equation}
\noindent where the $\mathbf{r}_{f_k}$ is the centre coordinate of the interface. Due to the property of the linear reconstruction, the value at
the vertex that is opposite of $f_3$ will be larger than the $q_{V(l)}^{\max}$ at the same vertex, which means the MLP condition is violated.

Similar deduction for the minimum constraint $q^{\min}(\mathbf{r})$ could be given, and thus the lemma is proved.  $\Box$

The original MLP condition was proved to be satisfying maximum/minimum principle, which guarantees the stability of MLP limiter \cite{Park2010}. By simply following
similar procedure, the maximum/minimum principle could be proved for weak-MLP condition.

\noindent \textbf{Theorem 1.}  For a following hyperbolic conservation law in two-dimensional space,
 \begin{equation} \label{eq:2D_Conser}
\frac{\partial q}{\partial t}+\frac{\partial f(q)}{\partial x}+\frac{\partial g(q)}{\partial y}=0,
\end{equation}
 \noindent finite volume method could be used for fully spatial-temporal discretization. While a monotone Lipschitz continue flux function is used for the
 calculation of numerical fluxes, a linear reconstruction that satisfies Eq.\ref{eq:wMLP} under an appropriate CFL condiction
  will satisfy the maximum/minimum condition.

\noindent \textbf{Proof.}  Assume the reconstruction is made on a triangular cell. The conservation law of Eq.\ref{eq:2D_Conser} could be discretized into a
semi-discrete form
 \begin{equation} \label{eq:sMLP_half}
|\Omega_i|\frac{\partial {q}_i}{\partial t}
  +\sum_{k=1}^{N_f}F(q_{k}^+,q_{k}^-)|\partial\Omega_{k}|=0.
\end{equation}

Then, the linear reconstruction satisfied the weak-MLP condition, $\nabla q_i$ and $\nabla q_j$, which reconstruct the interface values in the following form
\begin{equation}
q_{k}^+=q_i+\nabla {q}_i \cdot \Delta r_{ik}, \quad
{q}_{k}^-=q_j+\nabla {q}_j \cdot \Delta r_{jk},
\end{equation}
\noindent give the following inequality
\begin{equation}
\overline{q}_{k}^{\min} \le q_{k}^{\pm} \le \overline{q}_{k}^{\max}.
\end{equation}
 \noindent Then, by applying Eq.\ref{eq:f_extr}, following inequality is given
\begin{equation} \label{eq:ue0}
q^{\min}_{\mathcal{V}(i)}\le \overline{q}_{k}^{\min}   \le  {q}_{k}^\pm \le \overline{q}_{k}^{\max}\le q^{\max}_{ \mathcal{V}(i)}.
\end{equation}

Therefore, the proof of the \emph{Theorem} in subsection 3.2 of \cite{Park2012} could be followed, and then the final formula of maximum/minimum principle
\begin{equation}
 {q}_{ \mathcal{V}(i)}^{\min,n}\le q_i^{n+1} \le  {q}_{  \mathcal{V}(i)}^{\max,n},
\end{equation}
\noindent will be valid in CFL condition of
\begin{equation}
\Delta t\frac{L_i}{|\Omega_i|}\left(\sup\limits_{q_1,q_2 \in \left[ {q}_{ \mathcal{V}(i)}^{\min}, {q}_{\mathcal{V}(i)}^{\max}\right]}\left|\frac{\partial F}{\partial q_2}\left(q_1,q_2\right) \right|\right)\le\frac{1}{3},
\end{equation}
\noindent where $L_i$ is the perimeter of $\Omega_i$.
 $\Box$

The principle in three-dimensional circumstance could be further investigated in the framework of \cite{Park2010}. Therefore, by applying the weak-MLP condition, the numerical schemes will be less dissipative in computations.

\subsection{Strict MLP condition} \label{sec:strict}

\begin{figure}
\begin{center}
\subfigure[\label{fig:f:MLP_original}{Original condition}]{
\resizebox*{5cm}{!}{\includegraphics{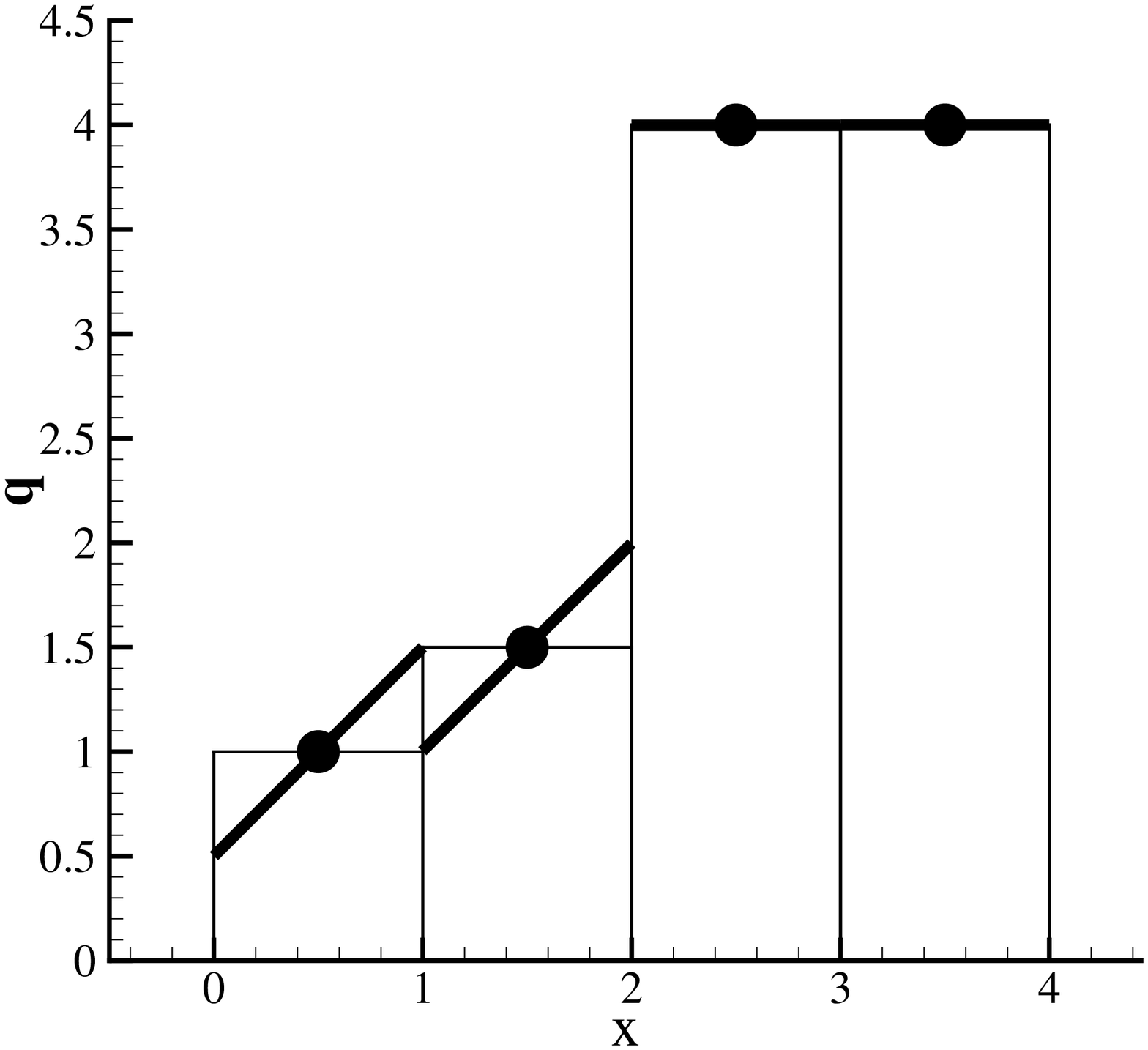}}} 
\subfigure[\label{fig:f:MLP_strict}{Strict condition}]{
\resizebox*{5cm}{!}{\includegraphics{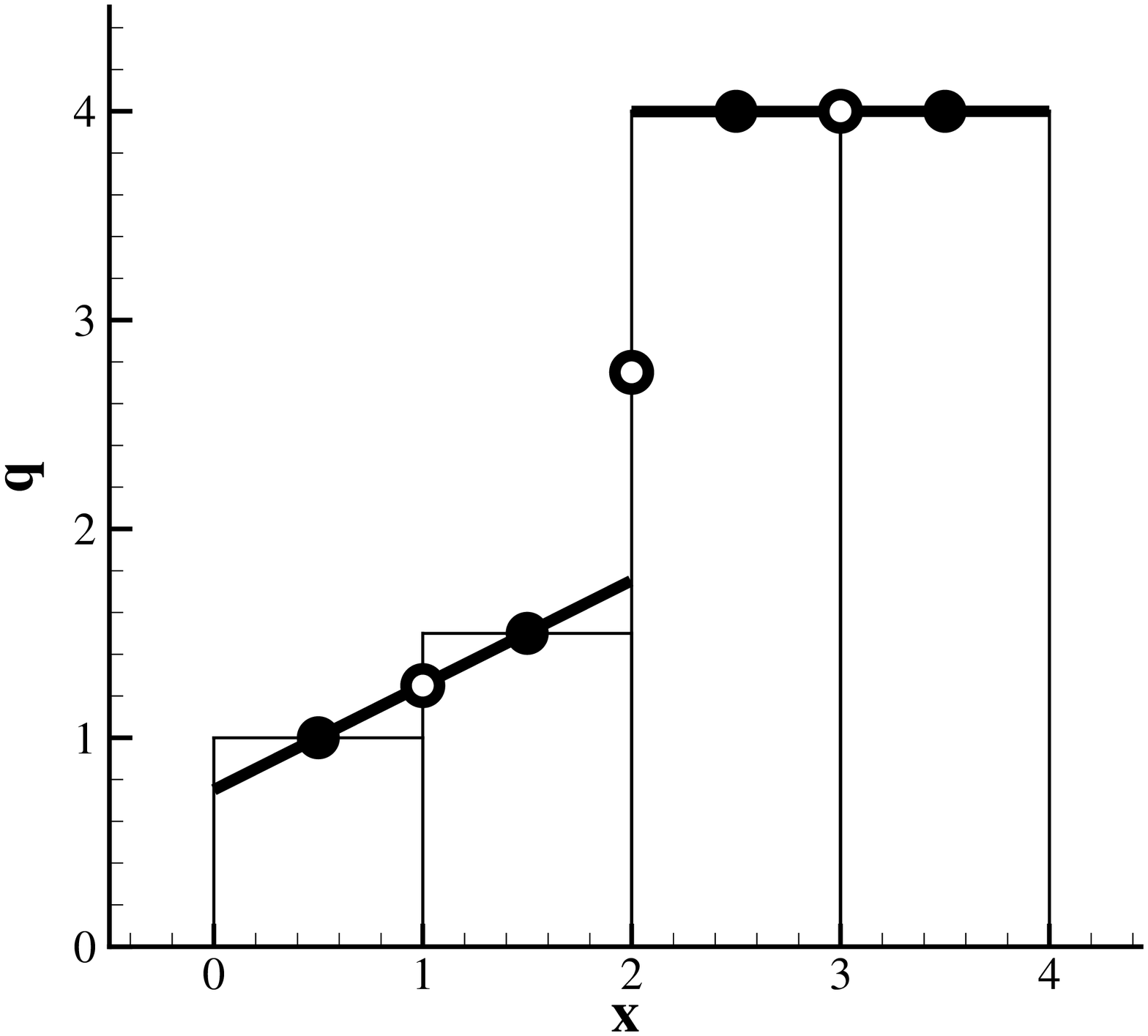}}}
\caption{\label{fig:f:Limiting} Limiting conditions based on different principles.}
\end{center}
\end{figure}

The strict MLP condition is defined for improving shock stability. By a direct observation, one may found out that the original MLP condition has a minor
 deficiency in certain situations, for which a example is shown in Fig.\ref{fig:f:Limiting}. In Fig.\ref{fig:f:MLP_original}, linear reconstructions are performed
 at each cell, and then limitations are made based on MLP condition. Therefore, new extrema (compared with cell centre values) will not be produced. In fact, in one-dimensional cases, Barth-Jespersen
 and Venkatakrishnan limiters present the same result.

In Fig.\ref{fig:f:MLP_original}, at the interface/point ($x=1$) between left two cells, the reconstructed value in the left side (of the interface/point) is larger than that in the right side. Whereas, the original values in cell centres are showing the contrary. Such a circumstance commonly exists in the conservative
piecewise-linear reconstruction of FVM which usually shows discontinuity at interfaces. Usually, this situation does not cause negative effects in computations.
However, it could cause instable results near shock waves because this unphysical distribution could affect the pre/pro-shock states.
As well known, shock wave is a highly non-linear phenomenon. A small change of pre/pro-shock states could change the strength or position of shock waves,
which could lead to remarkable flow variations.
Therefore, due to the nonlinearity, the potential effects should be aware of.

A strict monotonicity is hence presented in Fig.\ref{fig:f:MLP_strict}.
An average value is defined at each interface between two cells, and the average value will be used to define a bound for left or right reconstructed
values. For each cell, the limitation could be written as
\begin{equation} \label{eq:s_MLP}
\overline{q}^{\text{min}}_{v(i)}\le q \le \overline{q}^{\text{max}}_{v(i)}
\end{equation}
\noindent where the $\overline{q}^{\max,\min}_{v(i)}$ are the maximum/minimum of the averaged vertex values of cell $i$
\begin{equation} \label{eq:s_MLP_mm}
\overline{q}_{v(i)}^{\max}=\max\limits_{l\in v(i)}(\overline{q}_l), \quad
\overline{q}_{v(i)}^{\min}=\min\limits_{l\in v(i)}(\overline{q}_l).
\end{equation}
\noindent The overlines indicate that the variables are calculated by the average procedure, which could be a simple weighted average procedure as
\begin{equation}            \label{eq:Grad_NA}
\overline{q}_l=\frac{\sum_{i\in V(l)} \omega_{li} q_i}{\sum_{i\in V(l)} \omega_{li}}, \quad   \forall i, \omega_{li}\ge 0,
\end{equation}
\noindent where $\omega_{li}\ge 0$ is the weight of $q_i$. A simple and monotone weight is the inverse distance weight in \cite{Frink1991}, that is
\begin{equation}
\omega_{li}=1/|\Delta\boldsymbol{r}_{li}|.
 \end{equation}
 \noindent This average method will be used in the following paragraphs.

 Therefore, in any cases, the strict-MLP condition guarantees that the reconstructed values are strictly monotone. The maximum/minimum principle
is obviously satisfied by strict-MLP condition, which could be proved in a similar way as for weak-MLP condition. The brief proof is given as follow.

\noindent \textbf{Theorem 2.}   For a following hyperbolic conservation law in two-dimensional space,
 \begin{equation}
\frac{\partial q}{\partial t}+\frac{\partial f(q)}{\partial x}+\frac{\partial g(q)}{\partial y}=0,
\end{equation}
 \noindent finite volume method could be used for fully spatial-temporal discretization. While a monotone Lipschitz continue flux function is used for the
 calculation of numerical fluxes, a linear reconstruction that satisfies Eq.\ref{eq:s_MLP} under an appropriate CFL condiction
  will satisfy the maximum/minimum condition.

\noindent \textbf{Proof.}
Assume the reconstruction is made on a triangular cell. The semi-discrete conservation law in Eq.\ref{eq:sMLP_half} is also valid.

Then, the linear reconstruction satisfying the strict-MLP condition, $\nabla q_i$ and $\nabla q_j$, which reconstruct the interface values in the following form
\begin{equation}
q_{k}^+=q_i+\nabla {q}_i \cdot \Delta r_{ik}, \quad
{q}_{k}^-=q_j+\nabla {q}_j \cdot \Delta r_{jk},
\end{equation}
\noindent give the following inequalities
\begin{equation}
\overline{q}^{\min}_{v(i)} \le q_{k}^+ \le \overline{q}^{\max}_{v(i)}, \quad
\overline{q}^{\min}_{v(j)} \le {q}_{k}^- \le \overline{q}^{\max}_{v(j)}.
\end{equation}
 \noindent Then, by applying Eq.\ref{eq:f_extr}, following inequality is given
\begin{equation}
q^{\min}_{\mathcal{V}(i)}\le \overline{q}^{\min}_{v(i)}   \le  {q}_{k}^\pm \le  \overline{q}^{\max}_{ v(i)}\le q^{\max}_{ \mathcal{V}(i)},
\end{equation}
\noindent which is the same as that in Eq.\ref{eq:ue0}. Therefore, the subsequent conclusion could be made as that of Theorem 1.
$\Box$

\noindent \textbf{Remark 1.} The strict-MLP condition is more restrictive, and thus more dissipative. In multi-dimensional circumstances, it is unnecessary
to restrict the complete linear distribution within a cell to satisfy strict-MLP condition. In fact, only the values at several given point, for example the centres of
interfaces, should be checked for strict-MLP condition in Eq.\ref{eq:s_MLP}, and thus the scheme will be less dissipative.

\noindent \textbf{Remark 2.} Three conditions, original MLP, weak-MLP, and strict-MLP, which satisfy maximum/minimum principle, are already given.
By satisfying these three conditions, linear reconstructions could be expected to be free from spurious oscillations. However, as have already been proved, these conditions
are different in diffusivity, which will show different numerical performance.

\subsection{Pressure weight function} \label{sec:pwf}

Two new conditions are presented in the last two subsection, which will be used to reduce and enhance dissipation respectively. Therefore, how to combine these two conditions in a unified
framework is needed to be answered. Obviously, the strict-MLP condition is defined for improving shock stabilities, and thus this condition should be used in the vicinity
of shock waves. Flow pressure will be increasing drastically across shock waves, and thus the pressure increment had been used for indicating these strong discontinuities.
A effective pressure weight function is of the following form
\begin{equation} \label{eq:wf_MLP}
\omega_p= \left( \frac{\overline{p}^\text{min}_{v(i)}}{\overline{p}^\text{max}_{v(i)}}\right)^3,
\end{equation}
\noindent where the $\overline{p}^{\max,\min}_{v(i)}$ are maximum and minimum vertex pressure calculated by Eq.\ref{eq:s_MLP_mm} within a cell.

Similar polynomial type pressure weight function had been successfully applied for indicating shock waves in upwind schemes \cite{Kim1998,Zhang2016,Zhang2016a}.
This function will be a small value if the pressure variation within the cell is significant, and thus the strict-MLP condition could be applied.
On the contrary, weak-MLP condition will be used while the function gives a relatively large value (but never larger than one). Therefore, the following limiter is introduced
\begin{equation}
\begin{aligned}
\phi_\text{MLP-pw} =\min\limits_{ \{k|\partial \Omega_k \subset \partial \Omega_i\} } \left\{
\begin{aligned}
&f  \left(\frac{\omega_p\overline{q}_{k}^{\max}+(1-\omega_p)\overline{q}_{v(i)}^{\max}-q_i}{q_{k}-q_i} \right), &\text{if} \quad q_{k}-q_i > 0, \\
&f  \left(\frac{\omega_p\overline{q}_{k}^{\min}+(1-\omega_p)\overline{q}_{v(i)}^{\min}-q_i}{q_{k}-q_i} \right), &\text{if} \quad q_{k}-q_i < 0, \\
&1, &\text{if} \quad q_{k}-q_i = 0,
\end{aligned}
\right.
\end{aligned}
\end{equation}
\noindent where the $q_{k}$ is an interface centre value calculated by unlimited linear reconstruction from the centre of cell $i$.
To be specific, strict-MLP limiter can be attained by giving $\omega_p \equiv 0$, and weak-MLP limiter can be attained by giving $\omega_p \equiv 1$.

\noindent \textbf{Remark 3.} Here, pressure increment is utilised to indicate the existence of shock waves. It should be noted that density or entropy,
or velocity (vector) could be used for indicating discontinuities as well. For instance, in \cite{Kang2010}, density jump was used to distinguish linear discontinuities
from continuous region and nonlinear discontinuity. As well known, flows
across shock waves is showing entropy increasing.
Velocity vector is changed across shock waves, and thus it could be used to indicate shock wave, as in the rotated upwind scheme \cite{Ren2003,Nishikawa2008}.
Therefore, how to choose the shock indication parameter is an open question. Here, pressure is used because of its effectiveness. Furthermore,
the computations of contact discontinuity
or slip line are expected to be less dissipative, and thus the weak-MLP will be uniformly applied in continuous region and linear discontinuities.

 \noindent \textbf{Remark 4.} A polynomial function in Eq.\ref{eq:wf_MLP} is used to calculate the weight of strict/weak-MLP condition.
 The hyperbolic tangent function \cite{Kitamura2012a} and exponential function \cite{Wang2016} could serve similar purpose. However,
 the performance of these functions is not discussed in this article. The presented pressure weight function in Eq.\ref{eq:wf_MLP} will be showing
 satisfactory performance in the numerical cases.

  \noindent \textbf{Remark 5.} The accuracy of piecewise linear approximation is essentially second-order. The limiter methods will not elevate the accuracy but
  reduce it. By using the pressure weight function, the property of limiter will be varying depending on local flow field state. Obviously, the strict-MLP condition
  will produced more significant numerical dissipation, but the condition will only be activated within limited regions by the pressure weight function.

 \section{Numerical results} \label{sec:tests}

\subsection{Shock tube problems}

\begin{figure}
\begin{center}
\includegraphics[width=10cm]{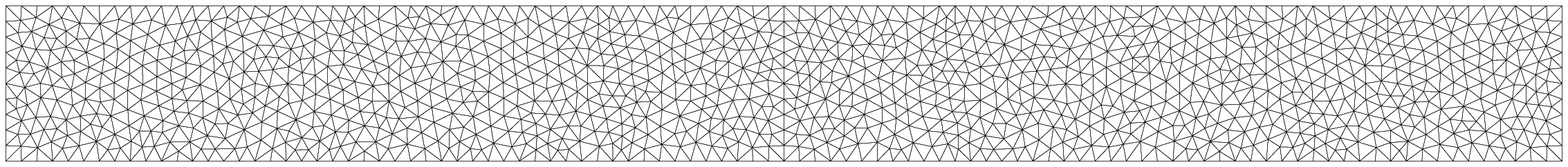}
\caption{Grid for the shock tube problems.}
\label{fig:f:ST_grid}
\end{center}
\end{figure}

Two shock tube problems are used to test limiters in simulations of some basic flow phenomenon on unstructured gird. The grid is shown in Fig.\ref{fig:f:ST_grid}. The computation domain is
$[0.0,1.0]\times[0.0,0.1]$. There are 101 boundary grid points in horizontal direction and 11 boundary grid points in vertical direction,
and 2292 triangular cells are created by Delaunay triangulation. A vertical grid line is formed in the middle of the domain, by
 which the initial discontinuity could be defined accurately. It should be noted that the grid is unsymmetrical, and thus the results will be unsymmetrical as well.
  Four steps Runge-Kutta scheme \cite{JAMESON1981,Frink1991} with CFL=0.2 is used for temporal solutions. Venkatakrishnan limiter,
  MLP limiter and MLP-pw limiter all use the Venkatakrishnan function $f_{\text{V}}$ used, of which the parameter $K$ is set as 1. HLLC scheme
\cite{Toro1994}  is used for the computations of convective fluxes.

\begin{figure}
\begin{center}
\subfigure[\label{fig:f:Sod_rho}{}]{
\resizebox*{5cm}{!}{\includegraphics{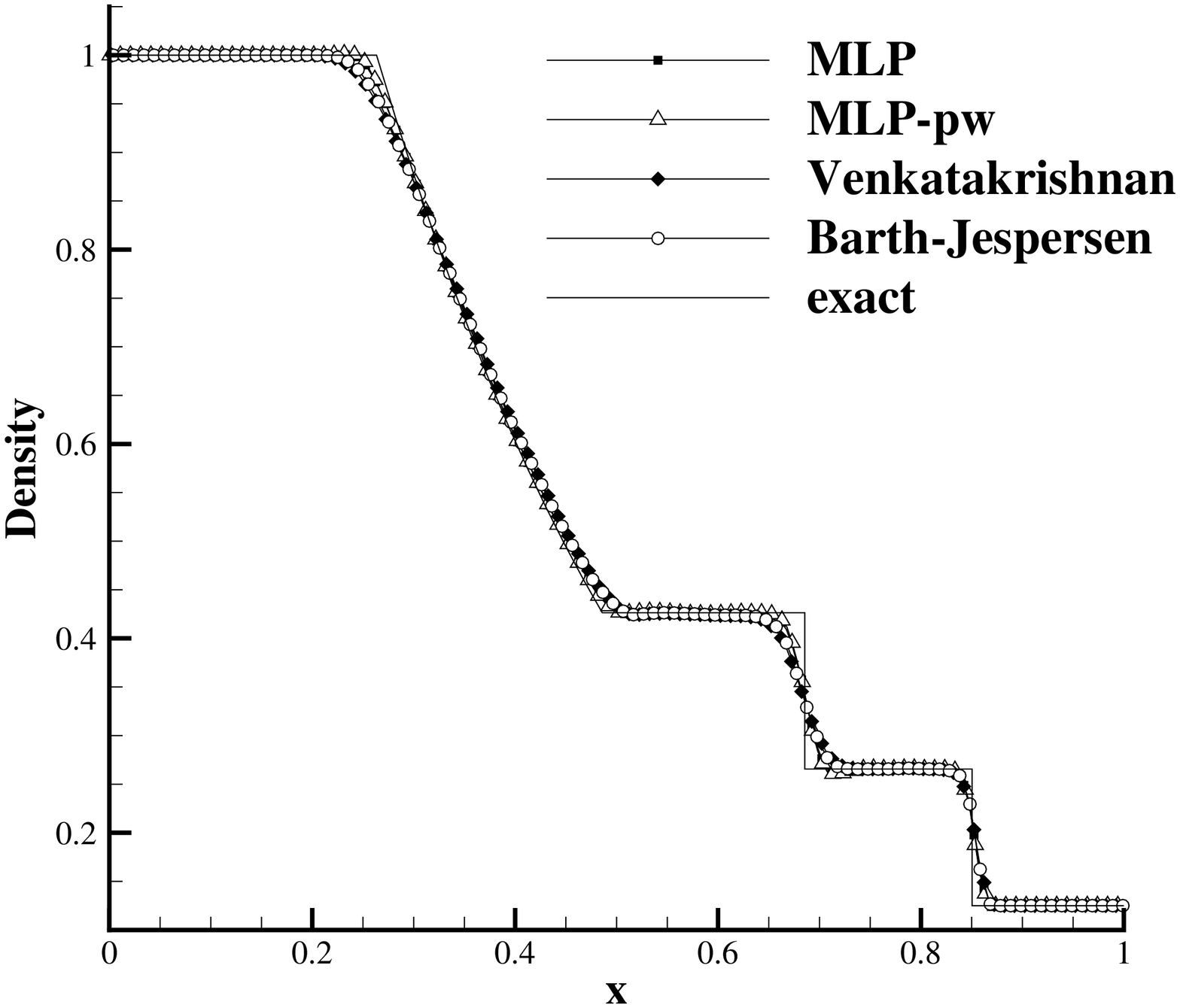}}} 
\subfigure[\label{fig:f:Sod_e}{}]{
\resizebox*{5cm}{!}{\includegraphics{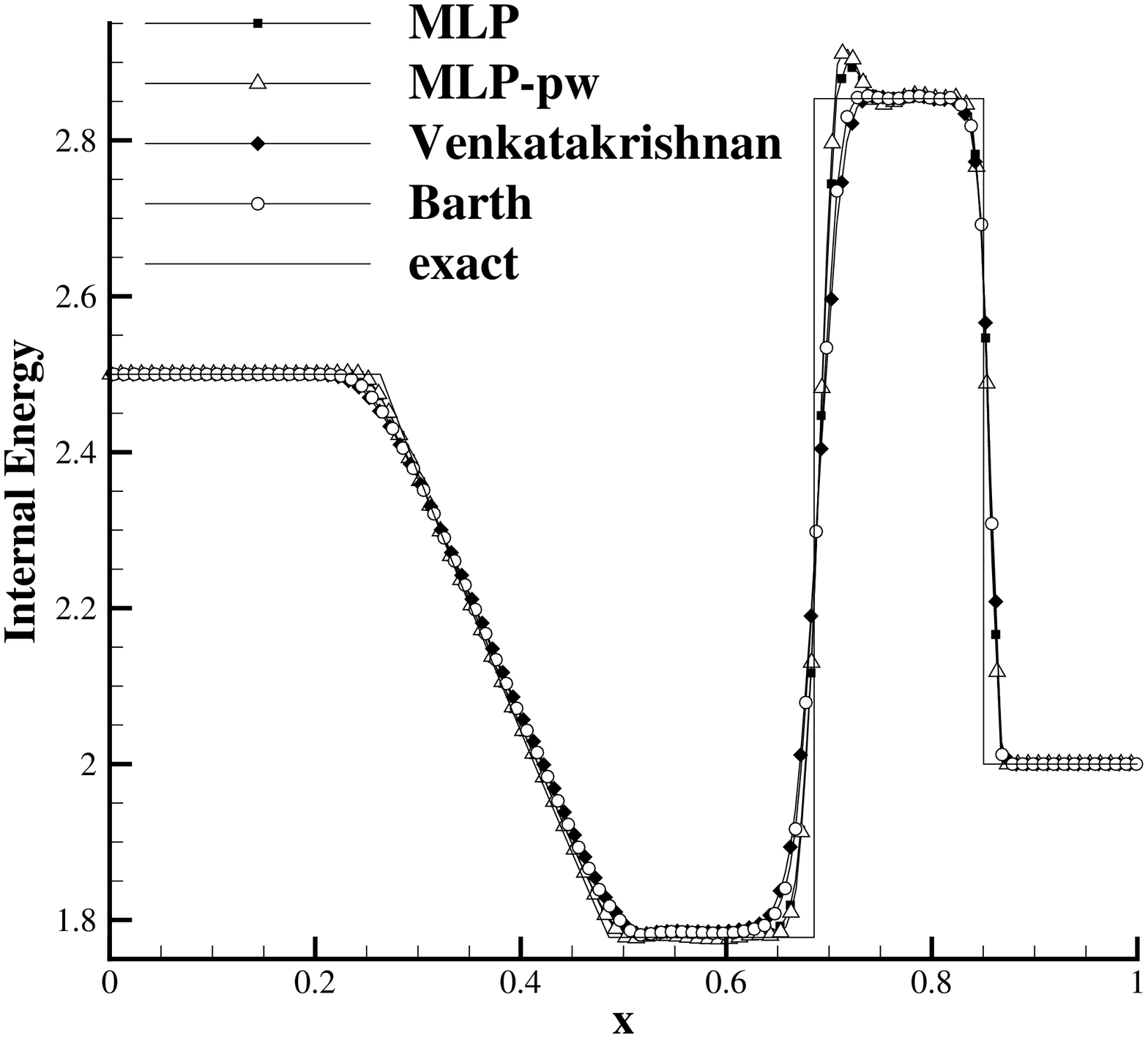}}}
\caption{\label{fig:f:Sod} Flow field distributions along the centerline of Sod problem.}
\end{center}
\end{figure}

The Sod shock tube problem \cite{Sod1978} is used for testing the performance of simulating shock wave, contact discontinuity and expansion.
The dimensionless initial condition across the middle discontinuity is $(\rho, u, v, p)_{L}=(1,0,0,1)$ and $(\rho, u, v, p)_{R}=(0.125,0,0,0.1)$. The results in non-dimensional time $t=0.2$ is presented in Fig.\ref{fig:f:Sod}.
Compared with Venkatakrishnan limiter or Barth-Jespersen limiter, MLP limiter and MLP-pw limiter are showing more accurate results,
but oscillations could be found in the right side of
contact discontinuity.

\begin{figure}
\begin{center}
\subfigure[\label{fig:f:Exp_rho}{}]{
\resizebox*{5cm}{!}{\includegraphics{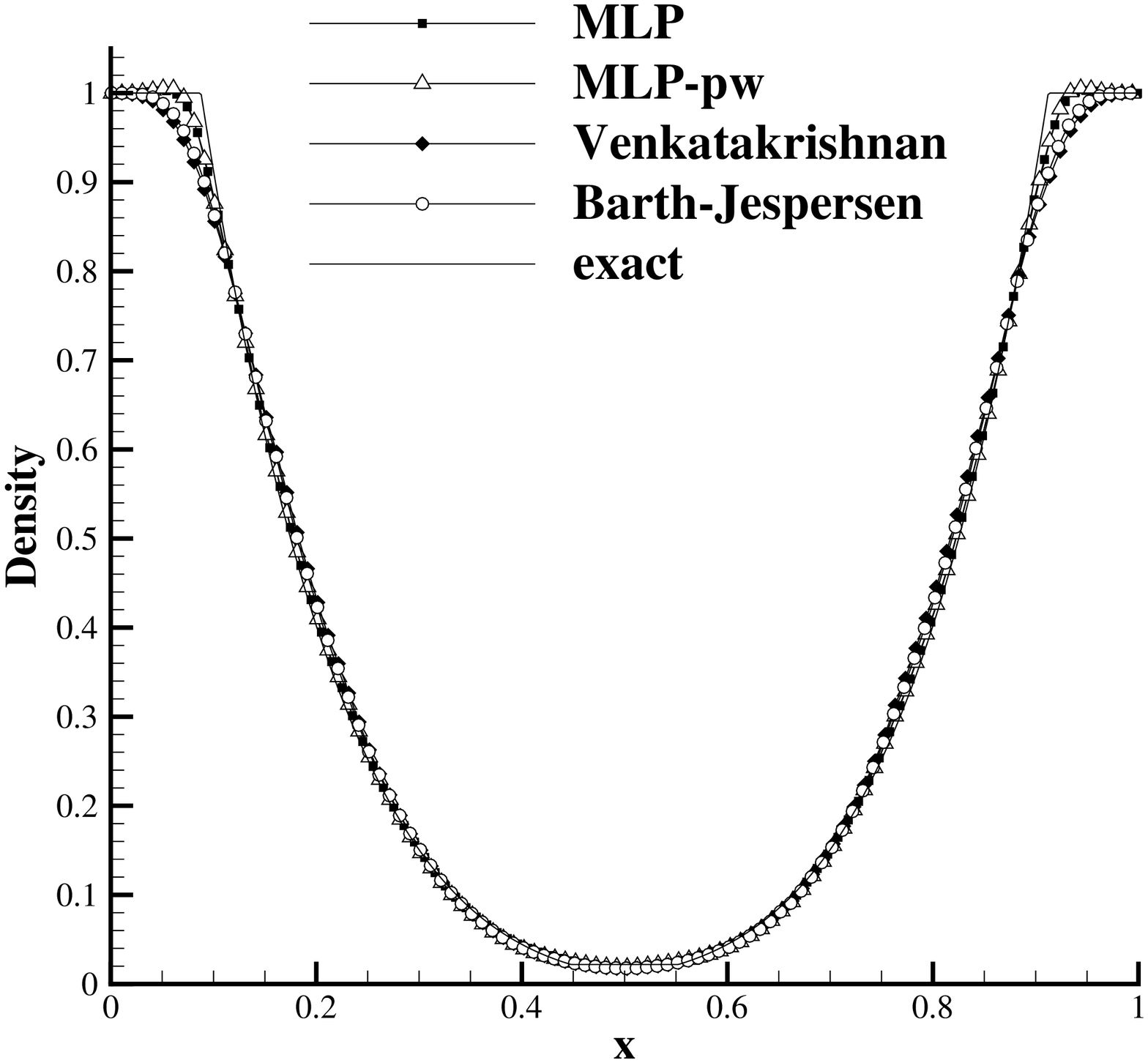}}} 
\subfigure[\label{fig:f:Exp_e}{}]{
\resizebox*{5cm}{!}{\includegraphics{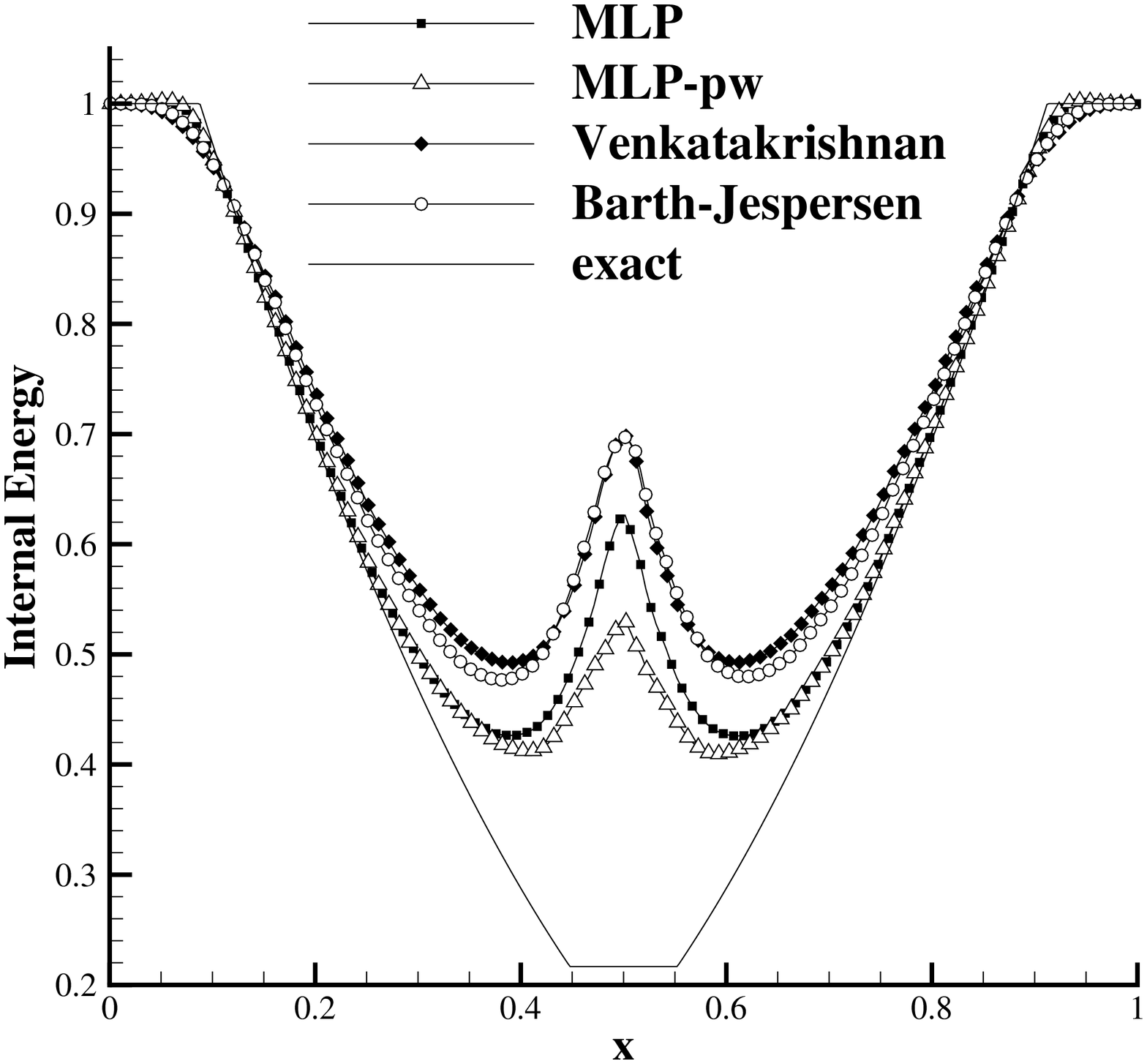}}}
\caption{\label{fig:f:Exp} Flow field distributions along the centerline of supersonic expansion problem.}
\end{center}
\end{figure}

Supersonic expansion problem is used to test the performance of computations in low density and pressure that approximate to zero.
The dimensionless initial condition across the middle discontinuity is $(\rho, u, v, p)_{L}=(1,-2,0,0.4)$ and $(\rho, u, v, p)_{R}=(1,2,0,0.4)$. The results in non-dimensional time $t=0.15$ is presented in Fig.\ref{fig:f:Exp}.
Again, Venkatakrishnan limiter and Barth-Jespersen limiter show more dissipative results.
MLP limiter and MLP-pw limiter are showing more accurate results, and the result calculated by MLP-pw limiter shows the most accurate approximation in the central low density region.

\begin{figure}
\begin{center}
\subfigure[\label{fig:f:sod_MLP_lim}{MLP}]{
\resizebox*{8cm}{!}{\includegraphics[angle=-90]{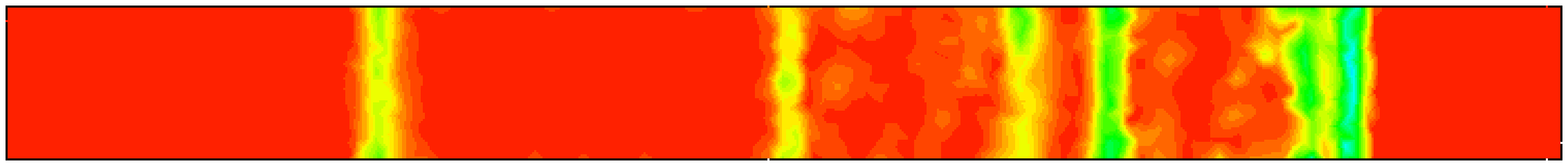}}} 
\subfigure[\label{fig:f:sod_MLPpw_lim}{MLP-pw}]{
\resizebox*{8cm}{!}{\includegraphics[angle=-90]{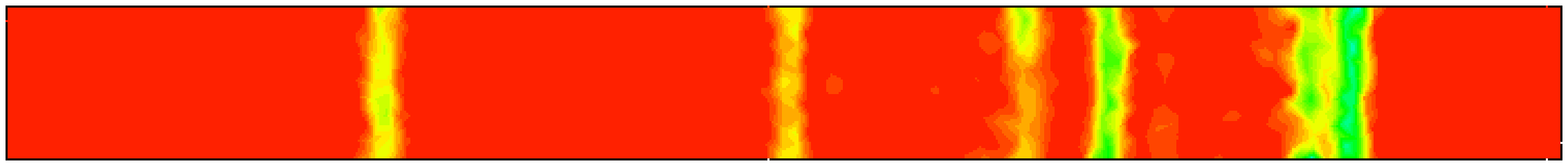}}}
\subfigure[\label{fig:f:sod_V_n_lim}{Venkatakrishnan}]{
\resizebox*{8cm}{!}{\includegraphics[angle=-90]{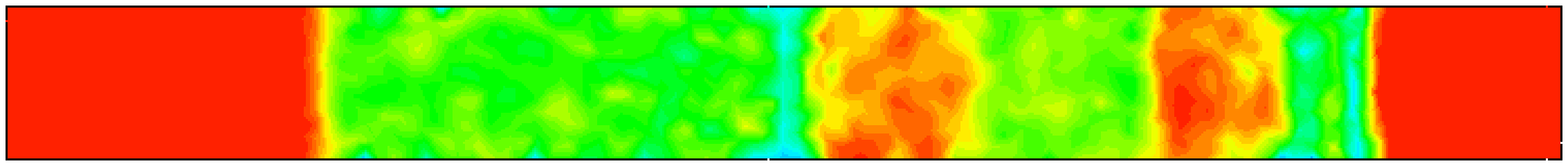}}} 
\subfigure[\label{fig:f:sod_B_n_lim}{Barth-Jespersen}]{
\resizebox*{8cm}{!}{\includegraphics[angle=-90]{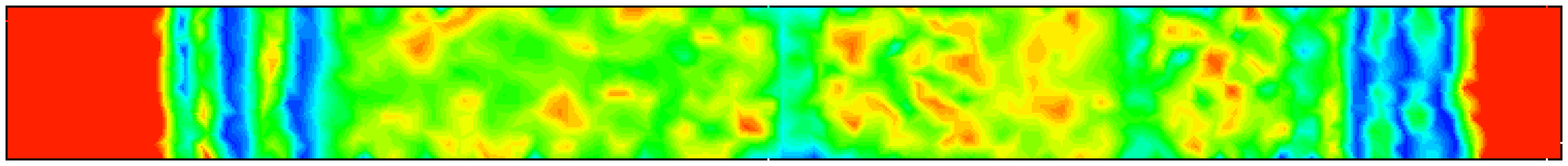}}}
\caption{\label{fig:f:sod_lim} Density limit value contours of Sod problem. Thirty equally spaced contour lines from $\phi=0$ (blue) to $\phi=1$ (red).}
\end{center}
\end{figure}

In order to reveal the detail of limitation process, density limit value contour of each limiter for Sod problem is shown in Fig.\ref{fig:f:sod_lim}. In general,
MLP-pw limiter is less diffusive in most of the regions, especially in continuous region and contact discontinuity. Correspondingly, oscillation has been found in contact discontinuity.
Similar behavior had occurred in the results of MLP limiter \cite{Park2010} and WBAP limiter \cite{Li2011}, for which the characteristic-wise reconstruction is
successfully used. Therefore, the oscillation of MLP-pw limiter has not been deem as a significant drawback, and more information will be given in the following results.

\subsection{Double shock reflection}

Inviscid supersonic flow over a wedge-shaped forward step is simulated, in which two reflection shock waves and an expansion fan are formed sequently. Flow field is discretized
by 58684 triangular cells. LU-SGS (Lower-Upper Symmetric-Gauss-Seidel) scheme \cite{YOON1988} is used for temporal solution and CFL=10. HLLC scheme is used for the computations of convective fluxes. The parameter $K$ of Venkatakrishnan function $f_{\text{V}}$ is set as 10. Mach number of uniform inflow is 2.

\begin{figure}
\begin{center}
\subfigure[\label{fig:f:DSR_MLP}{MLP}]{
\resizebox*{5cm}{!}{\includegraphics[angle=-90]{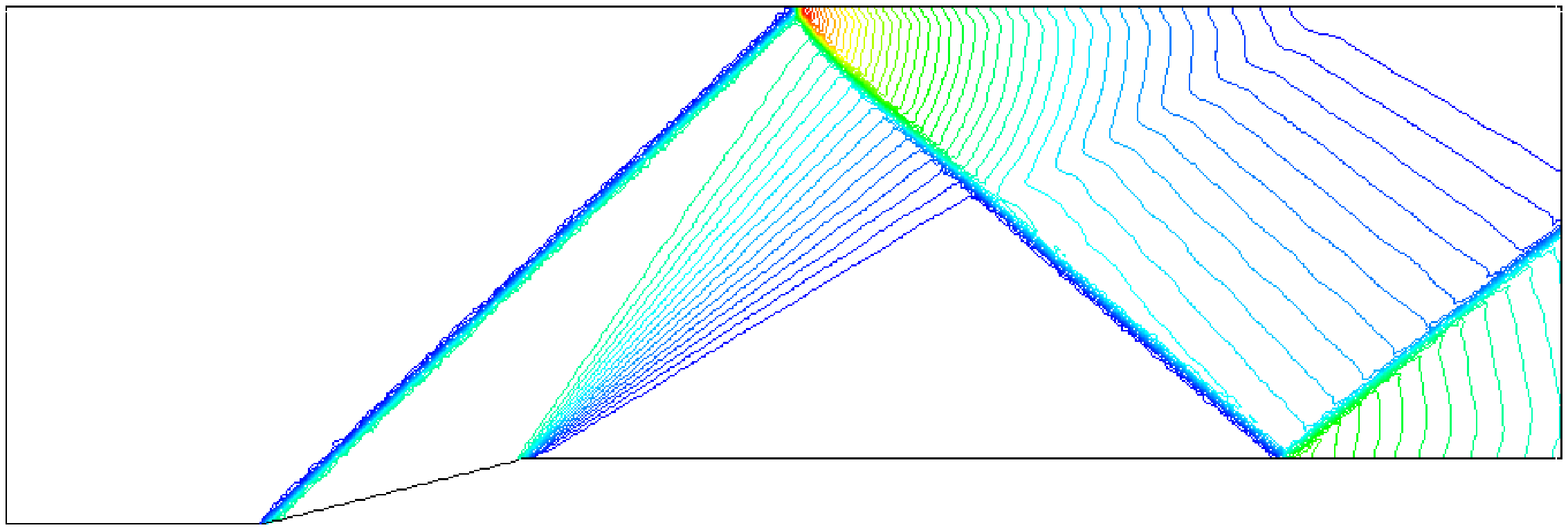}}} 
\subfigure[\label{fig:f:DSR_MLPpw}{MLP-pw}]{
\resizebox*{5cm}{!}{\includegraphics[angle=-90]{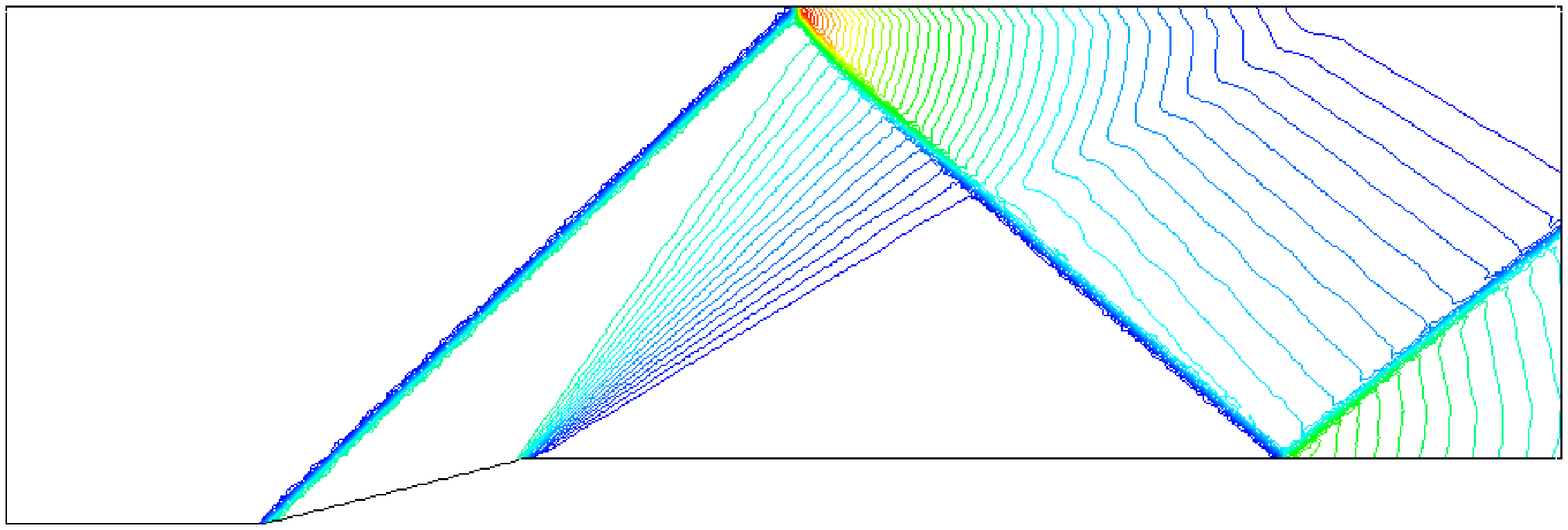}}}
\subfigure[\label{fig:f:DSR_V}{Venkatakrishnan}]{
\resizebox*{5cm}{!}{\includegraphics[angle=-90]{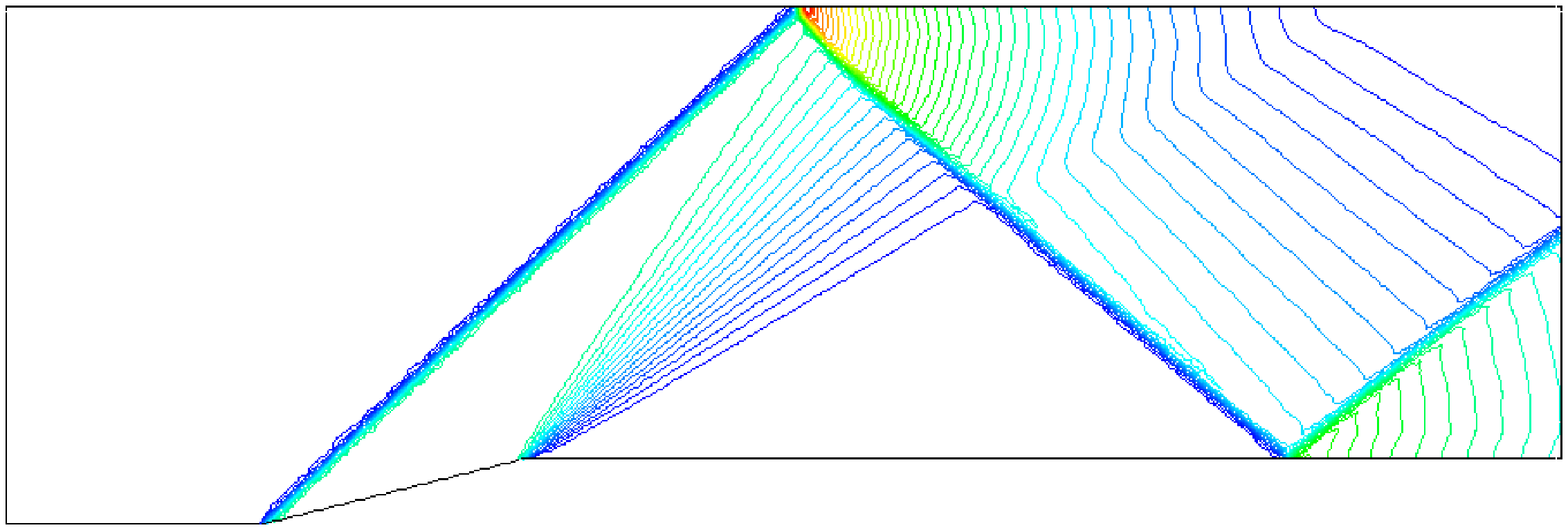}}} 
\subfigure[\label{fig:f:DSR_Barth}{Barth}]{
\resizebox*{5cm}{!}{\includegraphics[angle=-90]{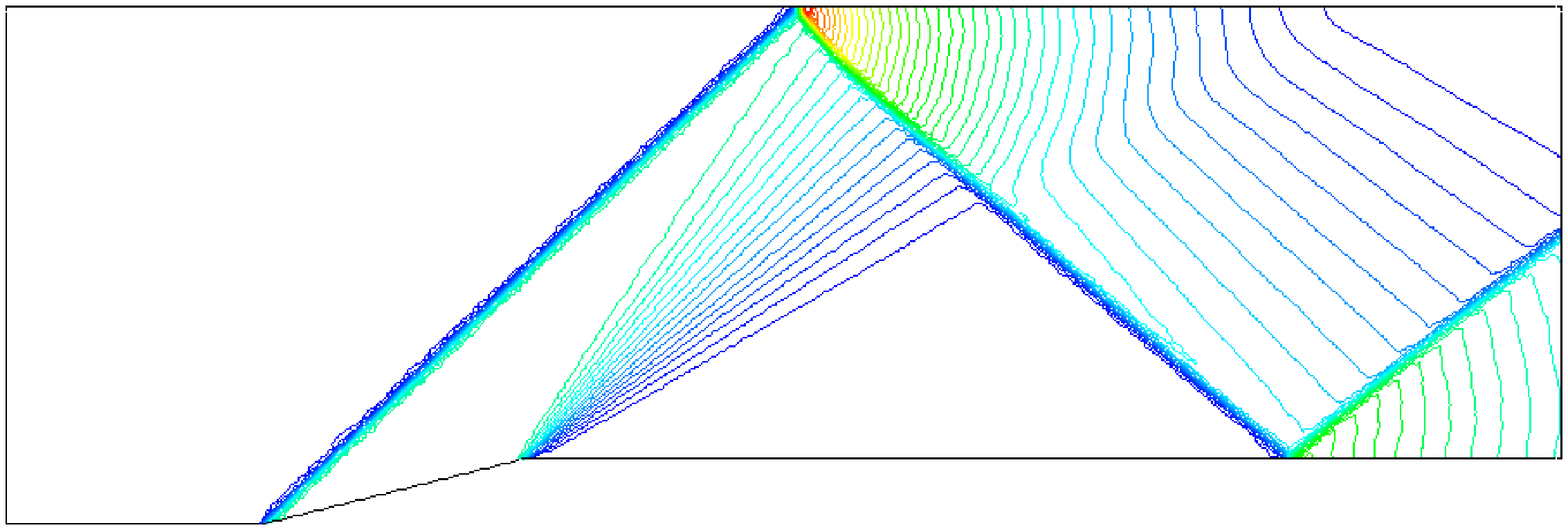}}}
\caption{\label{fig:f:DSR} Density contours of double shock reflection problem. Forty equally spaced contour lines from $\rho=1.0$ to $\rho=2.8$.}
\end{center}
\end{figure}

\begin{figure}
\begin{center}
\subfigure[\label{fig:f:DSR_MLP_lim}{MLP}]{
\resizebox*{5cm}{!}{\includegraphics[angle=-90]{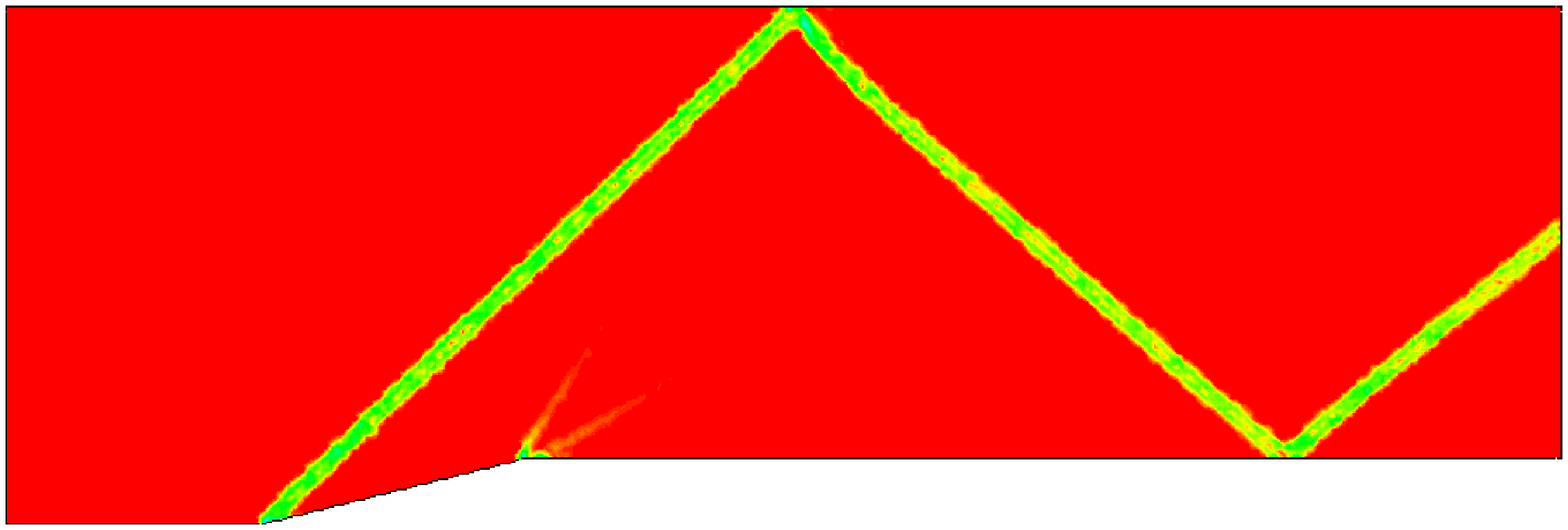}}} 
\subfigure[\label{fig:f:DSR_MLPpw_lim}{MLP-pw}]{
\resizebox*{5cm}{!}{\includegraphics[angle=-90]{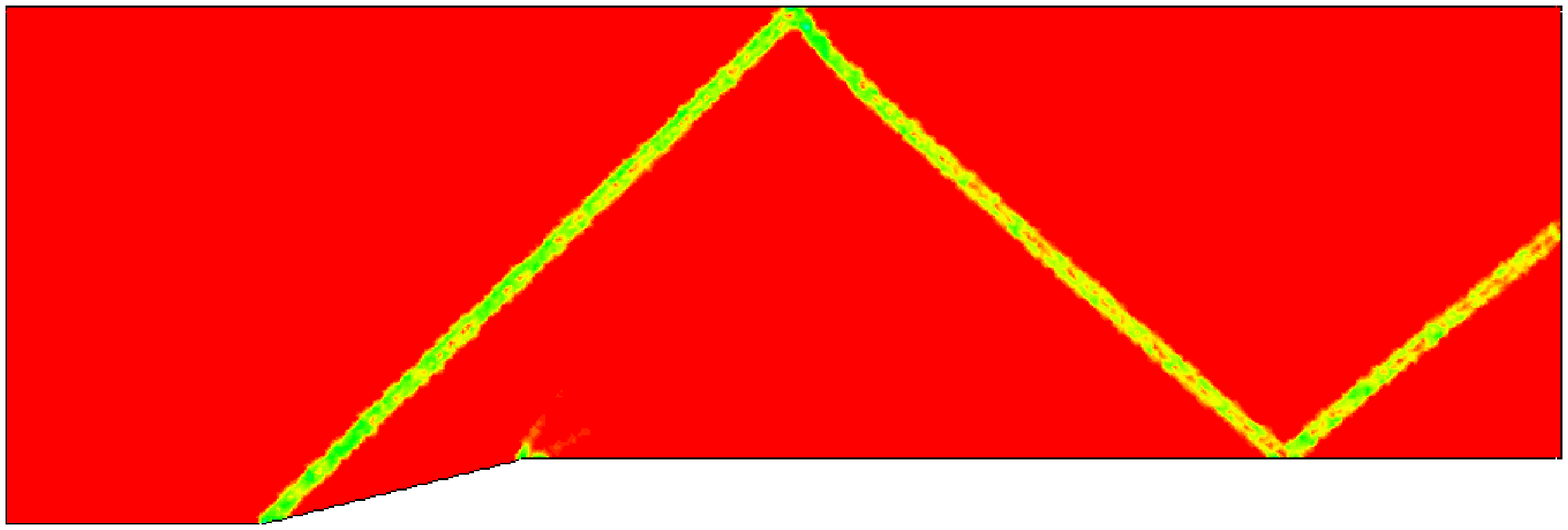}}}
\subfigure[\label{fig:f:DSR_V_lim}{Venkatakrishnan}]{
\resizebox*{5cm}{!}{\includegraphics[angle=-90]{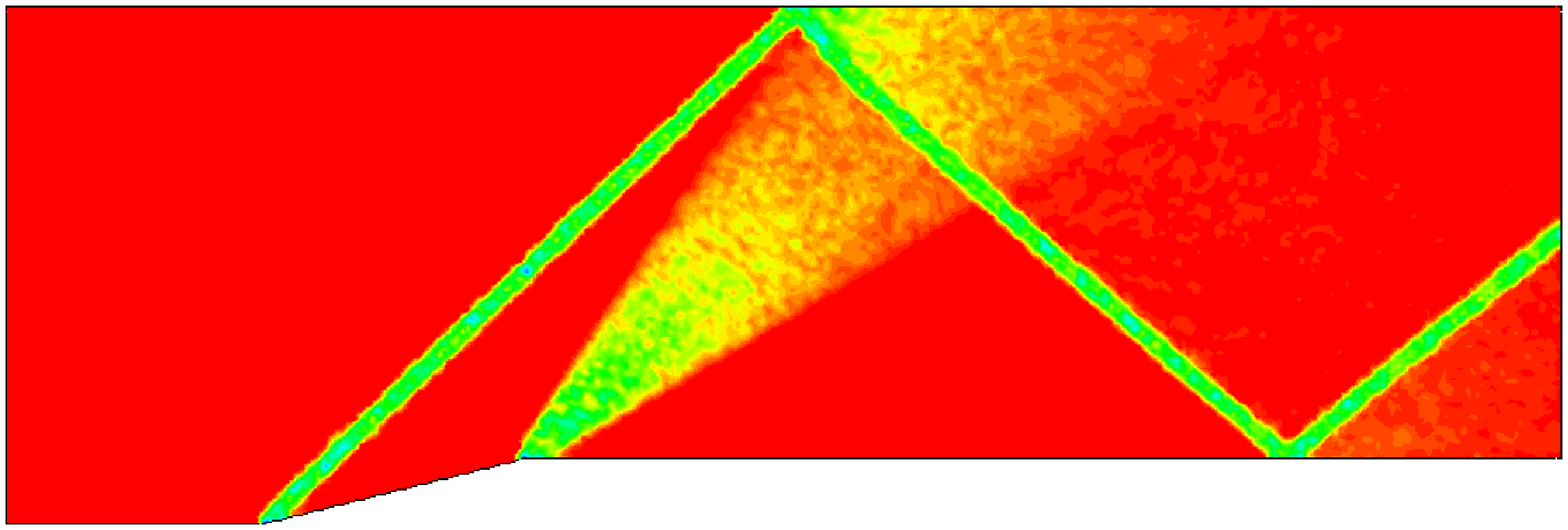}}} 
\subfigure[\label{fig:f:DSR_Barth_lim}{Barth}]{
\resizebox*{5cm}{!}{\includegraphics[angle=-90]{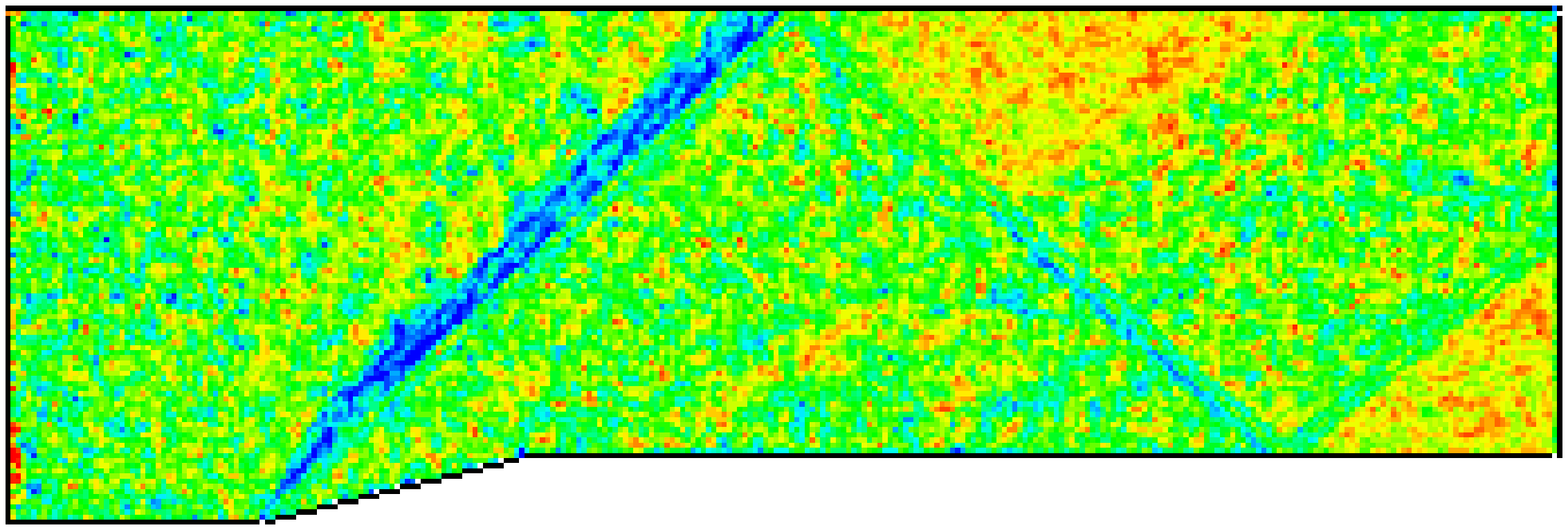}}}
\caption{\label{fig:f:DSR_lim} Density limit value contours of double shock reflection problem. Thirty equally spaced contour lines from $\phi=0$ (blue) to $\phi=1$ (red).}
\end{center}
\end{figure}

The density contours are shown in Fig.\ref{fig:f:DSR}. Slight differences could be found in the expansion fan.
The expansion fan in the results of Venkatakrishnan limiter or Barth-Jespersen limiter is relatively smeared.
The density limit value contours are shown in Fig.\ref{fig:f:DSR_lim}. It is obvious that Barth-Jespersen shows significant limitation in smooth region,
even in front of shock waves. Venkatakrishnan limiter narrows the limitation which is only found in shock waves and expansion fan. However,
it is unnecessary here to limit the gradient in expansion fan. MLP limiter and MLP-pw limiter show a little limitation in the corner causing expansion, and the limitation
 of MLP-pw limiter is less significant.

\begin{figure}
\begin{center}
\includegraphics[width=6cm]{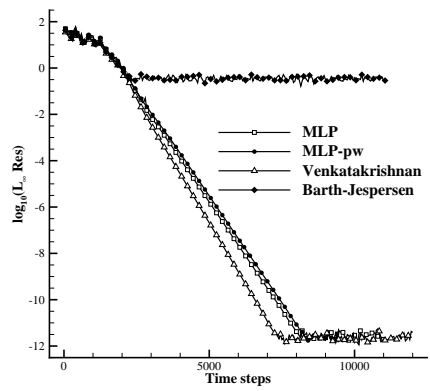}
\caption{Density residuals of the computations.}
\label{fig:f:DSR_res}
\end{center}
\end{figure}

The convergent histories are shown in Fig.\ref{fig:f:DSR_res}. Barth-Jespersen limiter is fail to be convergent.
Among the others, the Venkatakrishnan limiter shows slight advantage. The convergence performance of MLP and MLP-pw limiter is similar. In general, MLP-pw limiter
is less diffusive without depletion of convergence or stability.

\subsection{A Mach 3 wind tunnel with a step}
\begin{figure}
\begin{center}
\subfigure[\label{fig:f:Step_MLP}{MLP}]{
\resizebox*{5cm}{!}{\includegraphics[angle=-90]{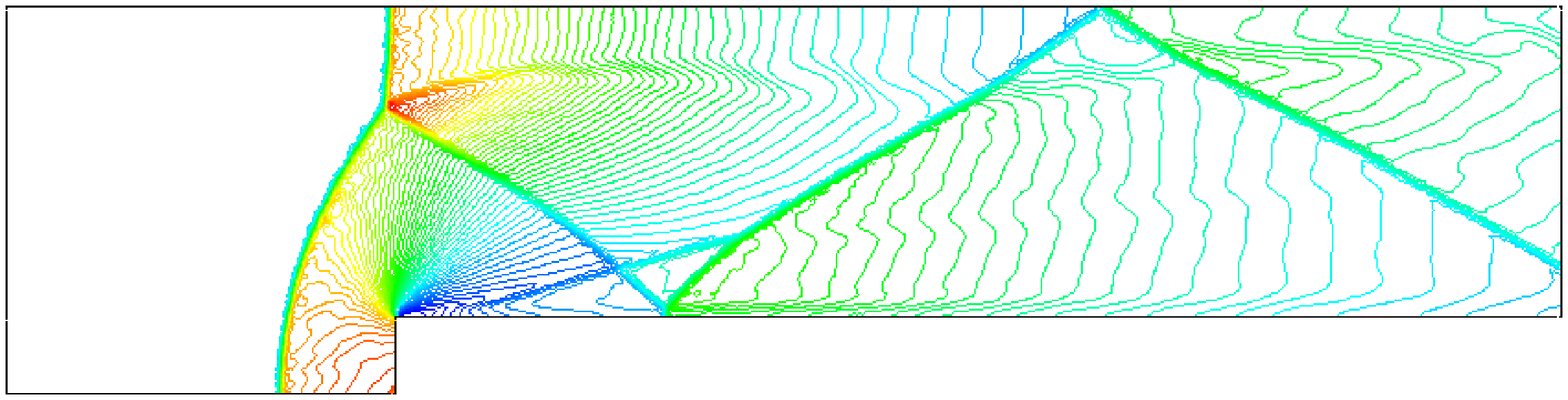}}} 
\subfigure[\label{fig:f:Step_MLPpw}{MLP-pw}]{
\resizebox*{5cm}{!}{\includegraphics[angle=-90]{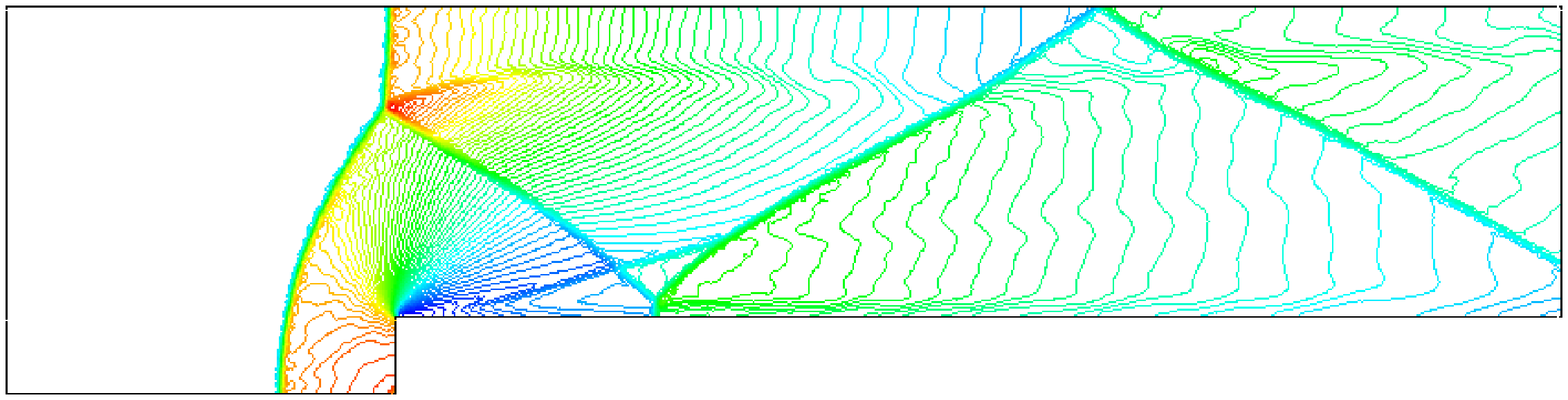}}}
\subfigure[\label{fig:f:Step_V}{Venkatakrishnan}]{
\resizebox*{5cm}{!}{\includegraphics[angle=-90]{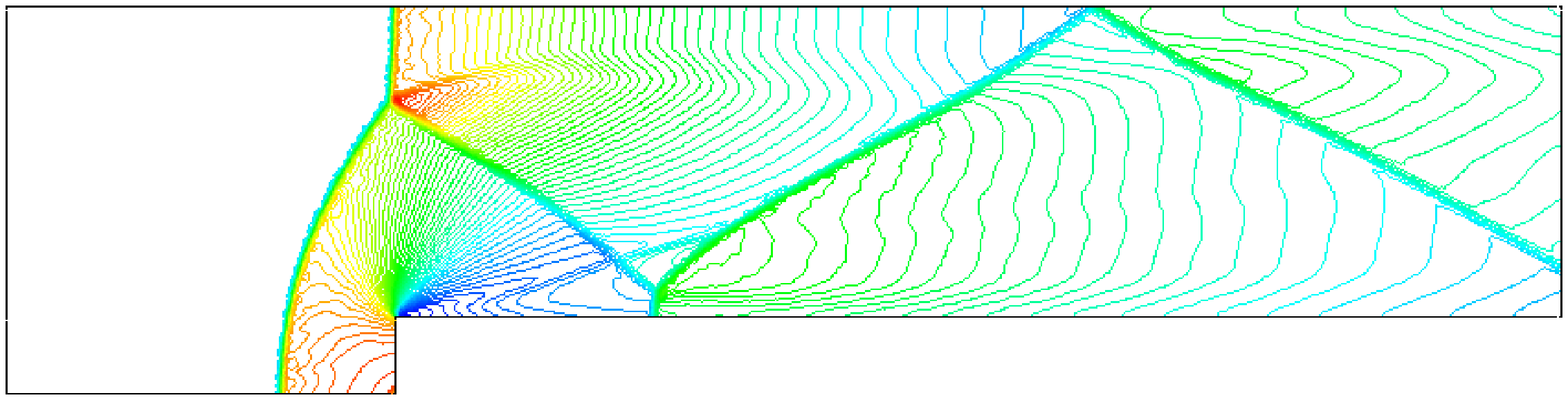}}} 
\subfigure[\label{fig:f:Step_Barth}{Barth-Jespersen}]{
\resizebox*{5cm}{!}{\includegraphics[angle=-90]{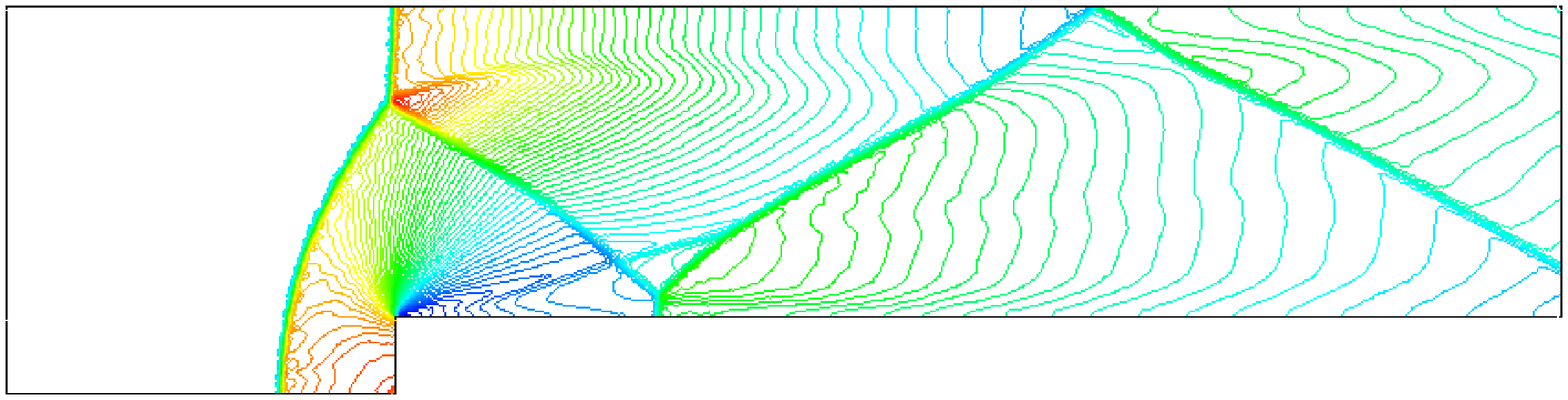}}}
\caption{\label{fig:f:Step} Density contours of the simulations of a Mach 3 wind tunnel with a step. Sixty equally spaced contour lines from $\rho=0.1$ to $\rho=4.5$.}
\end{center}
\end{figure}

\begin{figure}
\begin{center}
\subfigure[\label{fig:f:FS_s_MLP}{MLP}]{
\resizebox*{5cm}{!}{\includegraphics[angle=-90]{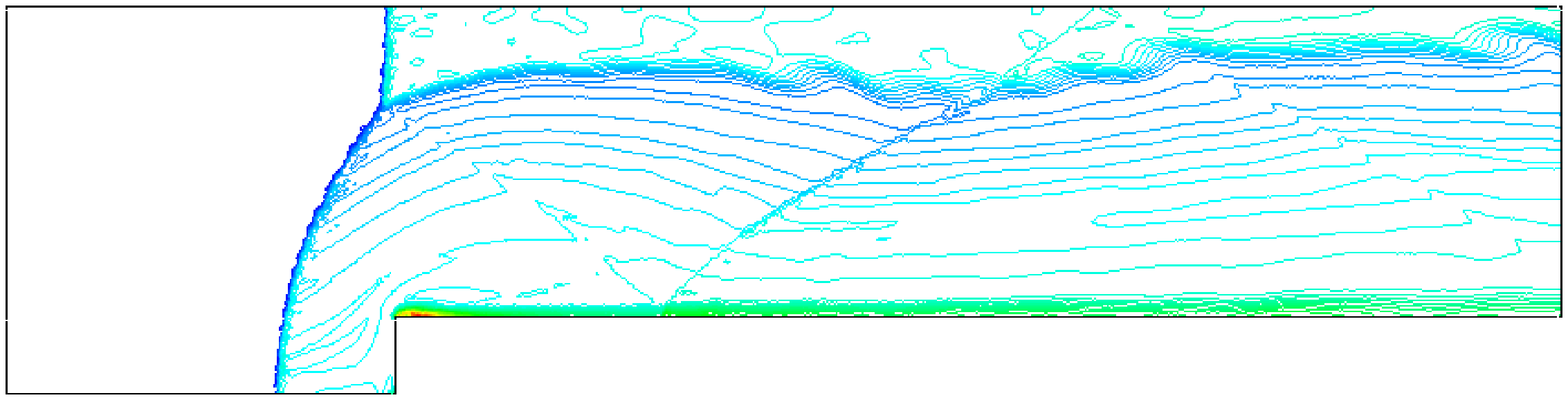}}} 
\subfigure[\label{fig:f:FS_s_MLPpw}{MLP-pw}]{
\resizebox*{5cm}{!}{\includegraphics[angle=-90]{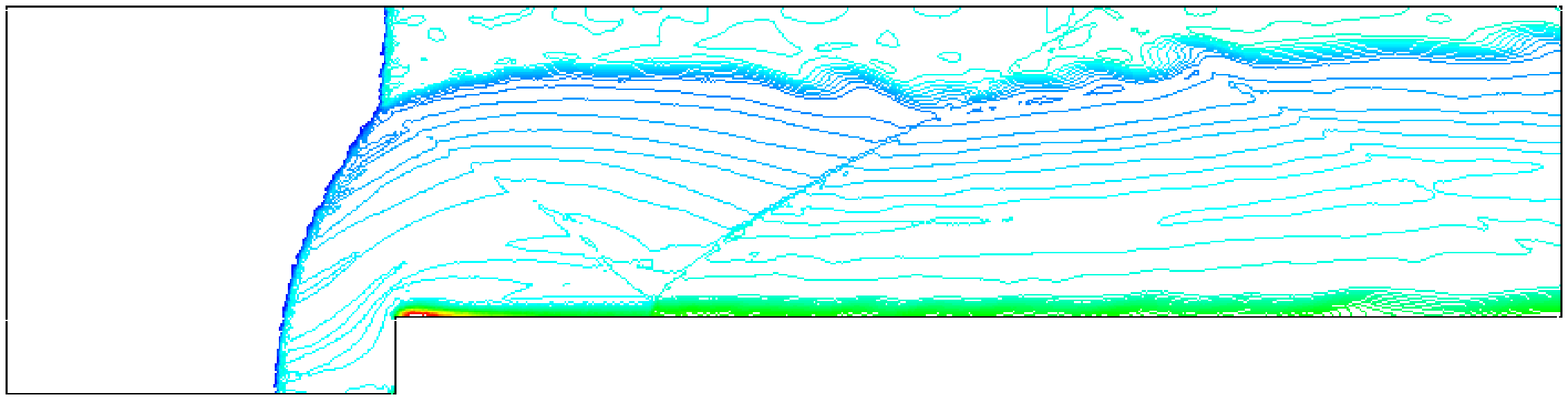}}}
\subfigure[\label{fig:f:FS_s_V_n}{Venkatakrishnan}]{
\resizebox*{5cm}{!}{\includegraphics[angle=-90]{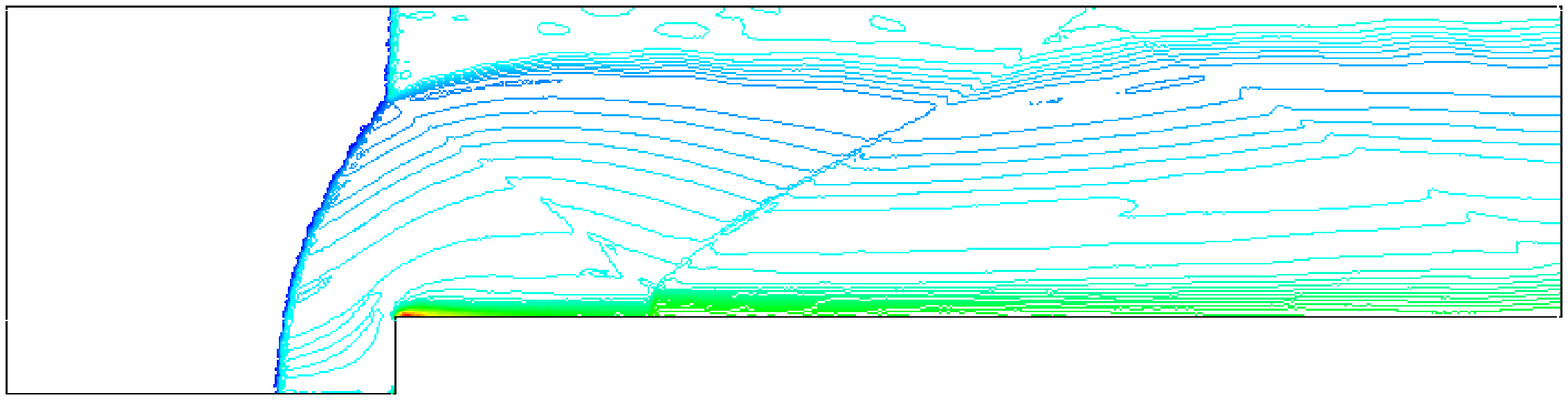}}} 
\subfigure[\label{fig:f:FS_s_B_n}{Barth-Jespersen}]{
\resizebox*{5cm}{!}{\includegraphics[angle=-90]{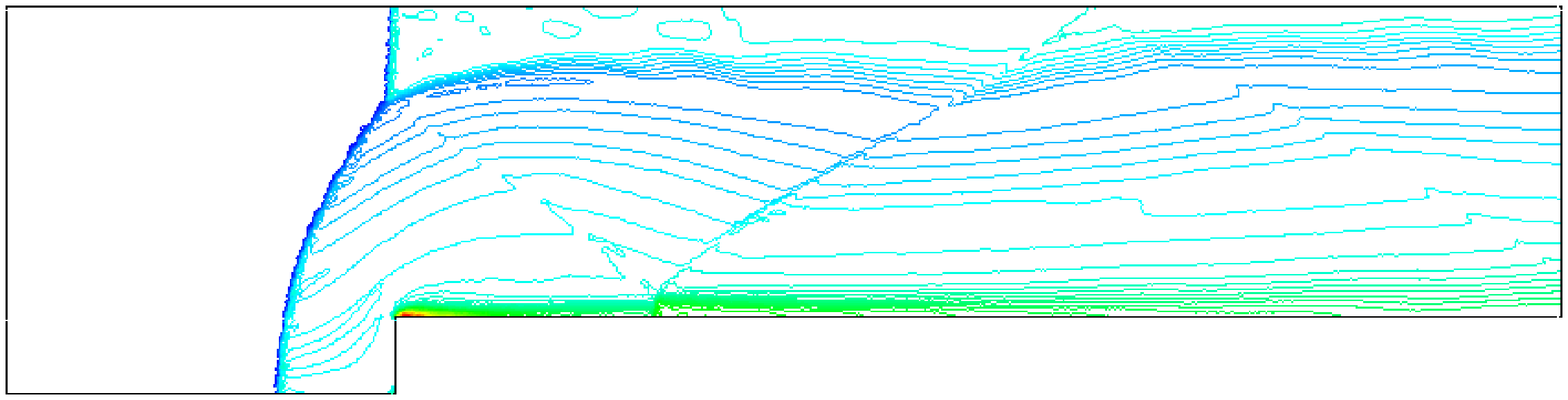}}}
\caption{\label{fig:f:FS_s} Entropy increment contours of the simulations of a Mach 3 wind tunnel with a step. Sixty equally spaced contour lines from $s=0.05$ to $s=2.05$.}
\end{center}
\end{figure}

Unsteady simulations are performed in this subsection. Uniform inviscid flow of which the Mach number is set as 3 passes a step will cause evolutive shock waves,
Mach stem and contact discontinuity \cite{Woodward1984}.
Flow field is discretized by 78246 triangular cells. Four-step Runge-Kutts scheme with CFL=1.5 is used for temporal solutions. AUSMPW scheme \cite{Kim1998} is used for
the computations of convective fluxes.
The parameter $K$ of Venkatakrishnan function $f_{\text{V}}$ is set as 10.
The results in non-dimensional time $t=4$ are presented in Fig.\ref{fig:f:Step}, and the entropy increment contours are shown in Fig.\ref{fig:f:FS_s}.
 The contact discontinuities of Venkatakrishnan limiter and Barth-Jespersen limiter are significantly smeared due to their higher dissipation.
 The contact discontinuity of MLP-pw limiter is more clear and develops to vortexes structure.
   Although not significant, the result of MLP limiter is more diffusive compared with that of MLP-pw limiter.

\begin{figure}
\begin{center}
\subfigure[\label{fig:f:Step_MLP_lim}{MLP}]{
\resizebox*{5cm}{!}{\includegraphics[angle=-90]{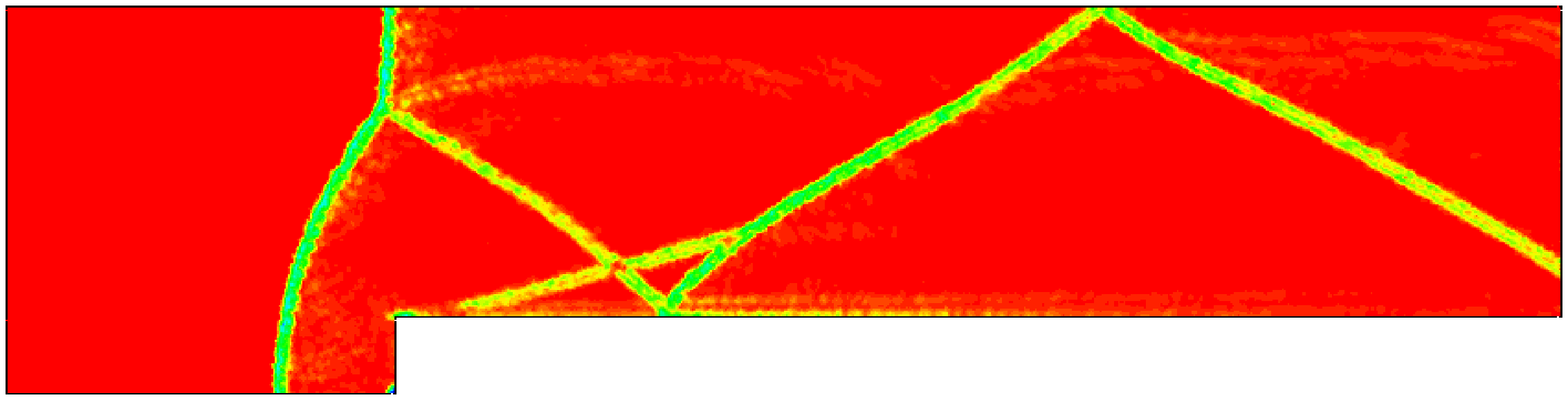}}} 
\subfigure[\label{fig:f:Step_MLPpw_lim}{MLP-pw}]{
\resizebox*{5cm}{!}{\includegraphics[angle=-90]{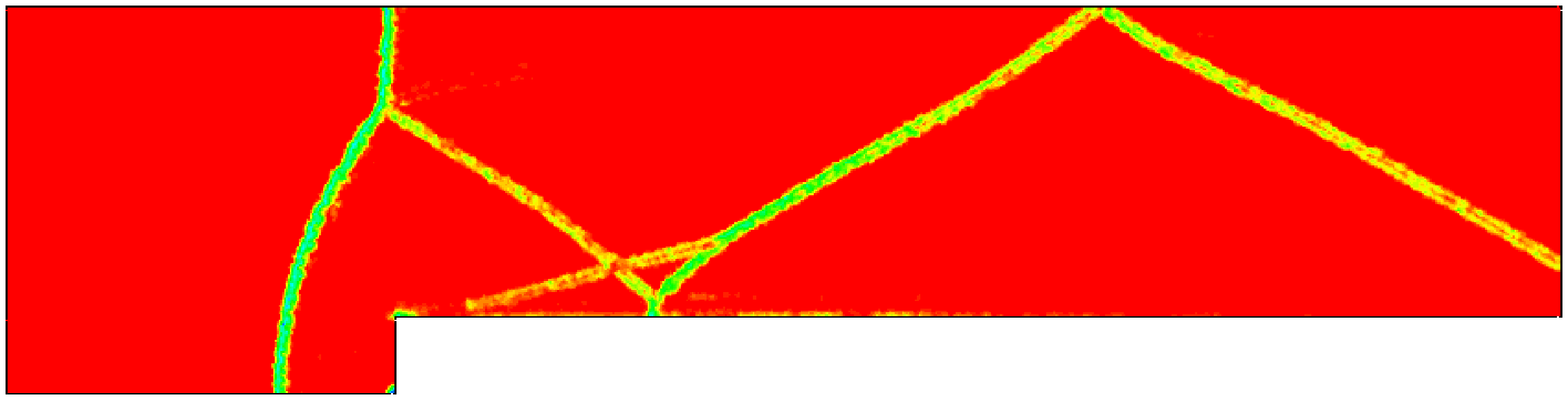}}}
\subfigure[\label{fig:f:Step_V_lim}{Venkatakrishnan}]{
\resizebox*{5cm}{!}{\includegraphics[angle=-90]{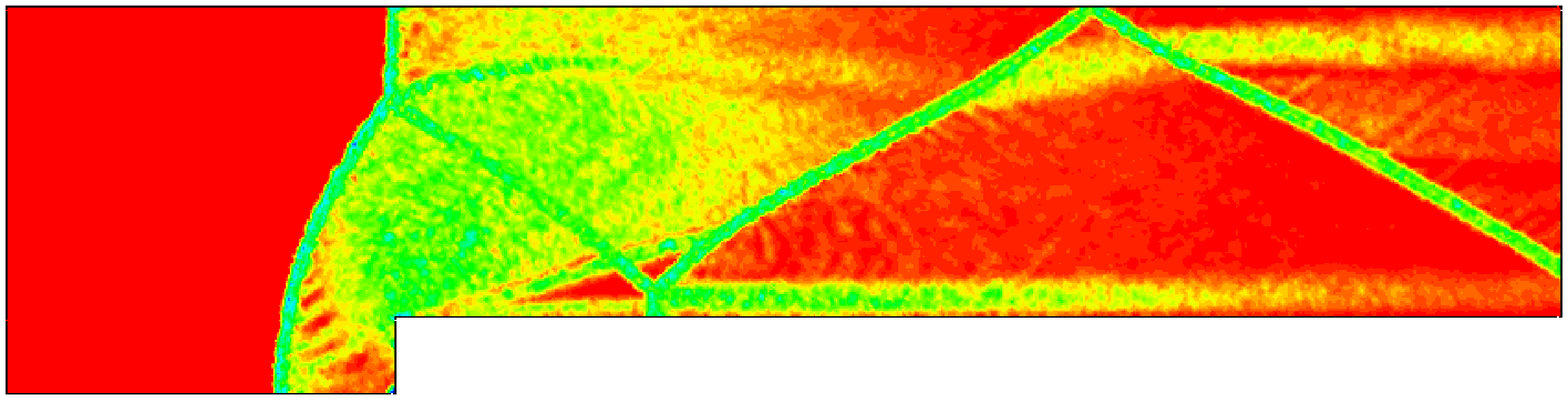}}} 
\subfigure[\label{fig:f:Step_Barth_lim}{Barth}]{
\resizebox*{5cm}{!}{\includegraphics[angle=-90]{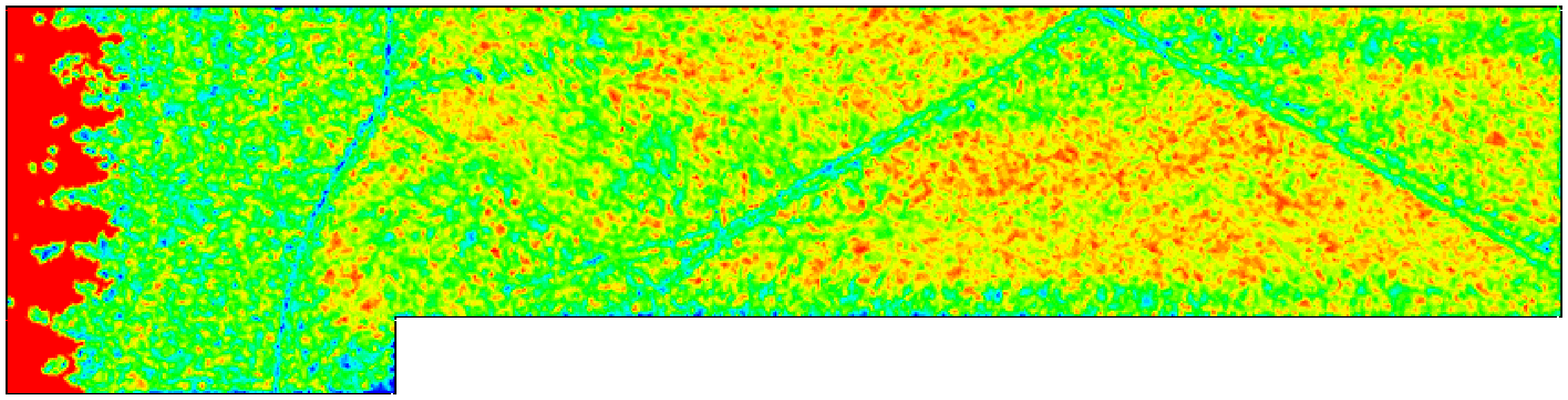}}}
\caption{\label{fig:f:Step_lim} Density limit value contours of the simulations of a Mach 3 wind tunnel with a step. Thirty equally spaced contour lines from $\phi=0$ (blue) to $\phi=1$ (red).}
\end{center}
\end{figure}

The density limit value contours are shown in Fig.\ref{fig:f:Step_lim}.
It could be found that Venkatakrishnan limiter and Barth-Jespersen limiter show extra dissipation on smooth region. MLP-pw limiter only shows limitation in
the vicinity of shock waves, which are pressure discontinuities. MLP limiter is less dissipative compared with Venkatakrishnan limiter and Barth-Jespersen limiter.
However, MLP still shows observable limitation in contact discontinuity and smooth region, and thus the density and entropy distributions are more diffusive.

\section{Conclusions} \label{sec:Conclusions}

As a modification of the multi-dimensional limiting process (MLP) on unstructured grids which significantly improves the convergence and accuracy of the simulations
on diverse problems, the presented method is a combination of two novel modified limit conditions, weak/strict-MLP condition.
Maximum/minimum principle are satisfied by both two new condition, and thus spurious oscillations will be well controlled. Especially,
the strict-MLP condition strictly limits the reconstructed variables, and thus the monotonicity could be guaranteed.
By using a pressure weight function that detects shock waves,
the strict-MLP condition is activated in the vicinity of shock waves and the weak-MLP condition is activated otherwise. Therefore, spurious oscillations are eliminated
near shock waves, even in hypersonic simulations, and the numerical dissipation are reduced in continuous regions and contact discontinuity.
Furthermore, the convergence of the presented limiter, MLP-pw, is improved.

\section*{Acknowledgments}
This work was supported by the National Natural Science Foundation of China under Grant 11372064, 91541117, and 11602052;
the Research Foundation of State Key Laboratory of Aerodynamics under Grant SKLA20160106; and
the National Key Research and Development Program under Grant 2016YFB0200700.
Xiao Liu at Dalhousie University gave advices to improve the writing of this article, who is sincerely appreciated.

\end{document}